\documentclass[11pt]{article}

\usepackage{amsmath,amsthm,amstext,amscd,amssymb,euscript,url}
\usepackage{mathrsfs}
\usepackage{calrsfs}
\usepackage{epsfig}
\usepackage[inline]{enumitem}
\usepackage[latin1]{inputenc}
\usepackage[T1]{fontenc}
\usepackage[english]{babel}
\usepackage{amsmath,amssymb,amsthm}
\usepackage{graphicx}
\usepackage{bbm}
\usepackage{color}
\usepackage{comment}
\usepackage[textsize=small]{todonotes}
\usepackage{color}
\usepackage{enumitem}

\newtheorem{theorem}{Theorem}[section]
\newtheorem{definition}{Definition}[section]
\newtheorem{lemma}{Lemma}[section]
\newtheorem{corollary}[theorem]{Corollary}
\newtheorem{proposition}[theorem]{Proposition}

\newtheorem{remark}[theorem]{Remark}

\newcommand{\be}{\begin{equation}}
\newcommand{\ee}{\end{equation}}
\newcommand{\nn}{\nonumber}

\newcommand{\Laplace}{\mathscr{L}}

\newcommand{\beq}{\begin{eqnarray}}
\newcommand{\eeq}{\end{eqnarray}}
\newcommand{\E}{\mathbb E}

\numberwithin{equation}{section}

\newcommand{\rem}[1]{}

\def\1{{\mathchoice {1\mskip-4mu\mathrm l}     {1\mskip-4mu\mathrm l}
{1\mskip-4.5mu\mathrm l} {1\mskip-5mu\mathrm l}}}
\newcommand{\indic}[1]{\1_{\{#1\}}}

\newcommand{\bl}{\begin{lemma}}
\newcommand{\el}{\end{lemma}}

\newcommand{\br}{\begin{remark}}
\newcommand{\er}{\end{remark}}

\newcommand{\bt}{\begin{theorem}}
\newcommand{\et}{\end{theorem}}

\newcommand{\bd}{\begin{definition}}
\newcommand{\ed}{\end{definition}}

\newcommand{\bind}{\begin{induction}}
\newcommand{\eind}{\end{induction}}

\newcommand{\bp}{\begin{proposition}}
\newcommand{\ep}{\end{proposition}}

\newcommand{\bc}{\begin{corollary}}
\newcommand{\ec}{\end{corollary}}

\newcommand{\bpr}{\begin{proof}}
\newcommand{\epr}{\end{proof}}

\newcommand{\bi}{\begin{itemize}}
\newcommand{\ei}{\end{itemize}}

\newcommand{\ben}{\begin{enumerate}}
\newcommand{\een}{\end{enumerate}}

\begin{document}






\title{
Large deviations and additivity principle \\for the open harmonic process}

\author{
Gioia Carinci$^{\textup{{\tiny(a)}}}$,
Chiara Franceschini$^{\textup{{\tiny(a)}}}$,
Rouven Frassek$^{\textup{{\tiny(a)}}}$,\\
Cristian Giardin\`a$^{\textup{{\tiny(a)}}}$,
Frank Redig$^{\textup{{\tiny(b)}}}$\;.
\vspace{.3cm}
\\
{\small $^{\textup{(a)}}$
FIM, University of Modena and Reggio Emilia},
{\small Modena, Italy}
\\
{\small $^{\textup{(b)}}$
Delft Institute of Applied Mathematics},
{\small TU Delft, Delft, The Netherlands}
}

\date{}
\maketitle

\pagenumbering{arabic}

\begin{abstract}
\noindent

We consider the boundary driven harmonic model, i.e. the Markov process associated
to the open integrable XXX chain with non-compact spins.
Using the factorial moments we characterize the
stationary measure as a mixture of product measures.
For all spin values, we identify the law of the mixture
in terms of the Dirichlet process. 
Next, by using the explicit knowledge of the
non-equilibrium steady state we establish formulas
predicted by Macroscopic Fluctuation Theory for several
quantities of interest: the pressure (by Varadhan's lemma),
the density large deviation function (by contraction principle),
the additivity principle (by using the Markov property of the mixing law).
To our knowledge, the results presented in this paper constitute the first rigorous
derivation of these macroscopic properties for models of energy transport
with unbounded state space, starting from the microscopic structure of the non-equilibrium steady state.
\end{abstract}

\newpage

\newpage

\section{Motivations and informal discussion of the main results}
In non-equilibrium statistical physics, a major  problem is to understand
systems with open boundaries, in particular the structure of their stationary
measure. In the literature this is often referred to as the ``non-equilibrium
steady state'' or the ``stationary non-equilibrium state''.
In the simplest set-up one considers one-dimensional models
on a finite segment of length $N$ which are driven out-of-equilibrium by two
boundary reservoirs with densities $\rho_l >0$, resp. $\rho_r>0$. 
A paradigmatic model,
for which  explicit knowledge of the stationary measure is available,
is the  boundary-driven simple symmetric exclusion process,  {where} one has
the description of the stationary measure via the matrix-product ansatz \cite{DEHP}.
Other models  are solvable but do not  {exhibit}
the long-range correlations structure that is believed to be a distinguishing feature
of non-equilibrium, such as zero-range  models
\cite{schutz, spitzer, zero-andjel}
which have a non-equilibrium steady state which is product,  or the Ginzburg-Landau model \cite{dop,cgp}, 
whose non-equilibrium steady state is a Gibbs measure
with  exponentially decaying correlations.
Clearly there is urgent need to identify other boundary-driven models
for which one has full control of the stationary state.
This is especially important to extract universal large scale properties
via the asymptotic analysis.

\vspace{.2cm}
In this paper we will prove that the family of boundary-driven model
 introduced in \cite{fgk} (called ``harmonic models'' because it involves
harmonic numbers) admits an explicit
description of the invariant measure for each system size $N$   {as a probabilistic mixture}.
This family of models, labelled by a parameter $s>0$,  emerged
 {as} the {\em integrable version} of the family of discrete Kipnis-Marchioro-Presutti
models \cite{kmp, Transport} (the two families share the same large scale behavior).
The root of the exact solvability of the harmonic models can be traced
back to the fact that they are related  to the open integrable XXX spin chain
with non-compact spins \cite{KRS,Lipatov, fadeev,Beisert,fgk}.
Remarkably, this spin chain is integrable for all spin values $s>0$
and thus the whole family of harmonic processes
is exactly solvable. See  \cite{fg} where  the moments and the stationary state were
obtained in closed-form. 

\vspace{.2cm}
Our first main result is presented in Theorem \ref{mixtmeasthm},
where we prove that the stationary measure of the
harmonic models is a ``mixture of inhomogeneous
Gibbs distributions''. A more precise, yet informal, version of this statement
is the following. Start from the equilibrium set-up (equal reservoir densities
$\rho_l=\rho_r$) and observe that the reversible Gibbs distribution
of the harmonic models is an homogenous product measure, the
marginal at each site being given by a Negative Binomial
distribution with shape parameter $2s>0$ and mean equal to the density of the
reservoirs.  Theorem \ref{mixtmeasthm} tell us that in a non-equilibrium set-up
(different reservoir densities $\rho_l\neq\rho_r$) the invariant measure of the harmonic models
is a mixture of inhomogeneous products of Negative Binomials
distributions with shape parameter $2s>0$ and {\em scale parameters which are given by
random variables, representing a random chemical potential at each site}.
We identify the law of these random
variables in terms of the symmetric Dirichlet distribution
with parameter $2s>0$ on the $(N+1)$-dimensional simplex.
As it is well know, when the parameter $2s$ is an integer, the Dirichlet distribution can
be expressed in terms of the order statistics of i.i.d. uniform random
variables. Our result agrees with the steady state obtained in  \cite{fg} (see Appendix
\ref{app:compare}), and reduces to the case of \cite{cfggt} where the stationary measure of the
harmonic model with s = 1/2 was proved to be a mixture of i.i.d.
geometric random variables whose mean are the order statistics of i.i.d.
uniforms.

\vspace{.2cm}
A second motivation of this paper is the Macroscopic Fluctuation Theory (MFT) \cite{mft},
which is a theory for diffusive systems proposed in recent years to describe the macroscopic
properties emerging in the limit $N\to\infty$.
MFT relies on the study of dynamical large deviations and states that macroscopically
the behavior of a diffusive systems is dictated by two transport coefficients,
the diffusivity $D(\rho)$ and the mobility $\sigma(\rho)$ depending on the system density $\rho:[0,1]\to\mathbb{R}_+$.
For the simple symmetric exclusion process, for which a dynamical
large deviation principle is available \cite{kov}, several findings of MFT
nicely match the results obtained with microscopic computations using
Bethe ansatz methods. See for instance \cite{DLS} for the large deviations
of the density profiles in the stationary state, \cite{DG1,DG2}
for the large deviations of the current and \cite{IMS} for the large deviations
of the positions of tagged particles. More recently, the time-dependent
solution of the MFT dynamical equations was found in \cite{MMS}
using  {integrability.}

\vspace{.2cm}
The boundary-driven harmonic models  considered in this paper,
labelled by a parameter $s>0$,  belongs to the class of models with
constant diffusivity and {convex} quadratic mobility
\be
\label{coeff}
D(\rho) = \frac{1}{2s} \qquad\qquad \text{and} \qquad\qquad \sigma(\rho) = \frac{\rho}{2s} \left(1+\frac{\rho}{2s}\right).
\ee
Other particle models in the same class 
include the symmetric inclusion processes
\cite{grv,Transport} and the discrete Kipnis-Marchioro-Presutti models \cite{kmp,Transport}.
For all these models, 
the state space is non-compact and the {\em dynamical large deviation
principle is not available}. 
The reason is that the stationary measures have exponential tails, and the proof of the dynamical large deviation principle, based on super-exponential replacement lemmas requires super-exponential tails of the stationary measures. This technical obstacle has so far not been overcome, and so all the results 
based on the Macroscopic Fluctuation Theory such as in \cite{bgl, bert} are conditional on the solution of this (highly non-trivial) technical issue.
This is also the case for the corresponding continuous models of energy transport,
namely the Kipnis-Marchioro-Presutti models \cite{kmp,Transport} (see also \cite{dfg} for recent results), the Brownian
energy processes \cite{gkrv} and the integrable heat conduction models
recently introduced in \cite{ffg}.  {We also mention \cite{cedric}, where a stochastic model of linear oscillators is studied and large deviations for the temperature profile in the non-equilibrium stationary state are analyzed.}
Therefore, for the class of models with
constant diffusivity and {convex} quadratic mobility
it is crucial to substantiate the predictions of MFT with microscopic computations,
which is the second aim of this paper.
In the rest of this introduction we give a summary of those MFT predictions,
first formulated in \cite{bgl} for the discrete Kipnis-Marchioro-Presutti model,
that we prove here for the boundary-driven harmonic model.

\vspace{.2cm}

\paragraph{Large deviations, pressure, additivity principle.}
We recall that, given a sequence of random variables $(X_n)_{n\ge 1}$
taking values in the measurable space $(\mathcal{X} ,\mathcal{B})$,
with $\mathcal{X}$ a topological space and $\mathcal{B}$ a $\sigma$-field of subsets of
$\mathcal{X}$, then we say that $(X_n)_{n\ge1}$ satisfies a large deviation principle with rate
function $I(x)$ and speed $N$ w.r.t. a sequence of probability measure $(\mu_n)_{n\ge 1}$ if,
for all $B\in\mathcal{B}$
$$
- \inf_{x\in B^o} I(x) \le
\liminf_{n\to\infty} \frac{1}{n} \log \mu_n(X_n \in B) \le
\limsup_{n\to\infty} \frac{1}{n} \log \mu_n(X_n \in B) \le
- \inf_{x\in \bar B} I(x)
$$
where $B^o$ denotes the interior of $B$ and $\bar B$ its closure.
Consider the empirical density profile
\be\label{edp}
L_N = \frac{1}{N}\sum_{i=1}^N \eta_i\, \delta_{\frac{i}{N}}
\ee
where $(\eta_i)_{i=1,\ldots,N}$ are distributed according to the invariant distribution
of the  boundary-driven harmonic model with parameter $s>0$,  system size $N\in \mathbb{N}$
and boundary densities $0\le\rho_l \le \rho_r <\infty$
(for a precise definition of the model see Section \ref{model-def}).
We introduce the space of density profiles
$$
\mathcal{X} = \{\rho \in L^1([0, 1], dx) : \rho(x) \ge 0\}
$$
equipped with the weak topology. Then, for models with transport coefficients \eqref{coeff} MFT predicts \cite{bgl} that
the sequence of empirical measures $(L_N)_{N\ge 1}$
satisfies a large deviation principle with speed $N$
and  rate function ${I}(\rho)$ which is the solution of the variational problem
\be
\label{rateI}
{I} (\rho) = \inf_{\theta} {\mathcal I} (\rho, \theta)
\ee
with
\be
\label{rateI2}
 {\cal I} (\rho, \theta) = 2s \int_{0}^1 dx \Big[
 \frac{\rho(x)}{2s} \log \frac{\frac{\rho(x)}{2s}}{\theta(x)}
 + (1+\frac{\rho(x)}{2s}) \log\Big(\frac{1+\theta(x)}{1+\frac{\rho(x)}{2s}}\Big)
 -  \log\Big(
 \frac{\theta'(x)}{\rho_r -\rho_l}\Big)\Big] \;.
\ee
The infimum in \eqref{rateI} is over increasing $C^1$ functions $\theta:[0,1]\to \mathbb{R}$  such that $\theta(0)=
\rho_l$ and $\theta(1)=  \rho_r$.
As remarked in \cite{bert} this large deviation function contains a relative entropy term and
a contribution related to the large deviations of the empirical profile of the order statistics of independent uniforms.
We will obtain rigorously (see  Theorem \ref{theo-ld}) this variational expression
from the exact description of the stationary measure, which indeed involves the order statistics
of independent uniforms.  
In particular the infimum in \eqref{rateI} corresponds to the contraction principle 
over the empirical profile of order statistics.

\vspace{0.1cm}
We will also study the pressure, which for a {function $h:[0,1]\to \mathbb{R}$} is defined as
\be
{P}(h) = \lim_{N\to\infty} \frac{1}{N} \log \mathbb{E} \Big[e^{N\langle L_N, h\rangle} \Big] \;.
\ee
The pressure can be obtained from the density large deviation rate function
via Legendre transformation, i.e., 
$$
{P}(h) = \sup_{\rho} \Big(  \int_{0}^1 h(x) \rho(x) dx  - {I}(\rho) \Big) \;.
$$
One gets the variational formula
\be
\label{pressureP}
{P}(h) =   \sup_{\theta}  {\cal P}(h,\theta)
\ee
with
\be
\label{pressureP2}
{\cal P}(h,\theta) =  2s  \int_0^1 dx \left[ \log \left(\frac{1}{1 + \theta(x)(1-e^{h(x)})}\right) +  \log\left(\frac{\theta'(x)}{\rho_r-\rho_l}\right)\right]
\ee
where again the supremum in \eqref{pressureP} is over increasing $C^1$ functions $\theta:[0,1]\to \mathbb{R}$
such that $\theta(0)= \rho_l$ and $\theta(1)= \rho_r$.
This will also be rigorously proved from the exact description of the stationary measure, see Theorem \ref{theo-pr}.
We remark that for models with constant diffusivity and {convex} quadratic mobility
it has been shown \cite{bgl} that the large deviation function of the density profile is non-convex and therefore the
Legendre transform of the pressure does not reproduce the large deviation function (it rather gives its convex hull).

\vspace{0.1cm}
Finally, the variational representations predicted by MFT encode an additivity principle \cite{bgl},
which can be formulated either for the pressure or for the density large deviation function.
For the pressure it is stated as follows.
For a macroscopic system of size $(b-a)$, where $-\infty<a<b<\infty$ with boundary parameters
$\rho_l,\rho_r$ define the modified pressure
$$
\widetilde P^{[a,b]}_{\rho_l,\rho_r}({h}) := P^{[a,b]}_{\rho_l,\rho_r}({h}) + 2s (b-a) \log\Big(\frac{\rho_r-\rho_l}{b-a}\Big)
$$
where
$$
P^{[a,b]}_{\rho_l,\rho_r}(h) = \sup_{\theta}
 \int_{a}^{b} dx
 \left[
2s \log \Big(
 \frac{1}{1+\theta(x) (1-e^{{h}(x)})} \Big)
+  \log\Big(
 \frac{(b-a)\theta'(x)}{\rho_r -\rho_{l}}\Big)
 \right] \;.
$$
Then, considering a macroscopic system of unit volume $[0,1]$ and
two subsystems of macroscopic  size $[0,x]$ and $[x,1]$ (with $0<x<1$),
the variational formula \eqref{pressureP}--\eqref{pressureP2}
of MFT is equivalent to the following additivity principle:
$$
\widetilde P^{[0,1]}_{\rho_l,\rho_r}({h})
=\sup_{\rho_l\leq \rho\leq \rho_r} \Big[\widetilde P^{[0,x]}_{\rho_l,\rho}({h}_1)  + \widetilde P^{[x,1]}_{\rho,\rho_r}({h}_2) \Big]
$$
where $h_1$ and $h_2$ are the restrictions of the function $h$ to the intervals $[0,x]$
and $[x,1]$.
Thus, the additivity principle relates the pressure of a macroscopic system of unit volume $[0,1]$ with boundary parameters
$\rho_l,\rho_r$ to the pressure of two subsystems, of macroscopic size $[0,x]$ and $[x,1]$
respectively, where the first subsystem is in contact with reservoirs of parameters
$\rho_l,\rho$ and the second subsystem  is in contact with reservoirs of parameters
$\rho,\rho_r$.
This will be proved in Theorem \ref{adprop} as a consequence of the
Markovian structure of the order statistics used to describe the stationary measure.
The additivity principle implies that the pressure of a constant field $h$, which corresponds to the large deviations of the total density, completely determines the pressure of any field. See 
 \cite{DLS} for the additivity principle of the density large deviation function of the symmetric exclusion process and \cite{bodineau} for a discussion of the additivity principle of the 
 time integrated current large deviation function, and its consequences in the setting of general diffusive systems.

\vspace{.3cm}
\noindent
{\bf Acknowledgements:}
In the course of this work we benefited from useful conversations with D. Gabrielli and D. Tsagkarogiannis, which were especially inspiring for the results in Section \ref{Mixt}.
C.G. would like to thank the Isaac Newton Institute for Mathematical Sciences, Cambridge, 
for support and hospitality during the SPL programme where work on this paper was undertaken;
in particular he acknowledges discussions with B. Derrida, T. Sasamoto and H. Touchette. 
This work was supported by Istituto Nazionale di Alta Matematica and EPSRC grant no EP/R014604/1. RF is supported in part by the INFN grant ``Gauge and String Theory (GAST)'' and by the ``INdAM-GNFM Project'', codice CUP-E53C22001930001.

\section{Model definition}
\label{model-def}

Denote by ${\Omega}_N$  the configuration space made of  $N$-dimensional vectors $\eta = (\eta_i)_{i\in\{1,\ldots,N\}}$
with non-negative integer components. We interpret the component $\eta_i$ as the number of particles
at site $i\in\{1,\ldots,N\}$. We shall write $\delta^i \in \Omega_N$ for the vector with all components zero
except in the $i^{th}$
place, i.e.
\be
({\delta}^i)_j=
 \left\{\begin{array}{rl}
 1 & \text{if } j=i,\\
 0 & \text{otherwise}.
 \end{array} \right.
\ee

\begin{definition}[Boundary-driven harmonic process with parameter $s>0$,  \cite{fg}]
\label{process}
For  $N\in\mathbb{N}$, we define the open symmetric harmonic process with parameter $s>0$
and  reservoir densities $0\le\rho_l \le \rho_r < \infty$ as the continuous-time Markov chain
$\{\eta(t)\,,\,t\ge 0\}$ having configuration space $\Omega_N$
and whose time-evolution is defined by the generator ${\mathscr L}$
working on functions $f : {\Omega}_N\to\mathbb{R}$
\be
\label{gen}
{\mathscr L}f :=  {\mathscr L}_1f + \Big(\sum_{i=1}^N {\mathscr L}_{i,i+1}f\Big)+ {\mathscr L}_Nf
\ee
where
\begin{eqnarray}
\label{gen-bulk}
({\mathscr L}_{i,i+1} f)(\eta)
& := &
\sum_{k=1}^{\eta_i}\varphi_s(k,\eta_i) \Big[f(\eta-k \delta^i + k \delta^{i+1}) - f(\eta)\Big] \\
& + &
\sum_{k=1}^{\eta_{i+1}} \varphi_s(k,\eta_{i+1}) \Big[f(\eta+k\delta^i - k \delta^{i+1}) - f(\eta)\Big] \nonumber
\end{eqnarray}
and, for $i\in\{1,N\}$,
\begin{eqnarray}
\label{gen-reservoir}
({\mathscr L}_{i} f)(\eta)
& := &
\sum_{k=1}^{\eta_i}\varphi_s(k,\eta_i) \Big[f(\eta-k \delta^i) - f(\eta)\Big] \\
& + &
\sum_{k=1}^{\infty} \frac{1}{k} \left(\frac{\rho_i}{1+\rho_i}\right)^k  \Big[f(\eta+k\delta^i) - f(\eta)\Big] \nonumber
\end{eqnarray}
with $\rho_1 = \rho_l$ and $\rho_N=\rho_r$.
Here the function $\varphi_s : \mathbb{N} \times \mathbb{N} \to \mathbb{R}$ is given by
\be
\label{varphi}
\varphi_s(k,n) :=  \frac{1}{k}\frac{\Gamma (n+1) \Gamma (n-k+2 s)}{ \Gamma (n-k+1) \Gamma (n+2 s)} \indic{1\le k\le n} \;.
\ee
\end{definition}

\medskip

\begin{remark}[Harmonic numbers]
{\em
When the occupation of the $i^{th}$ site is $n$, the function $\varphi_s(k,n)$ in \eqref{varphi} represents the rate at
which $k$ particles (with $1\le k \le n$)  jump from site $i$ to a nearest neighbour site $i\pm 1$.
One can  check that
\be
\label{eq:identity}
 \sum_{k=1}^n  \varphi_s(k,n) = \sum_{k=1}^n\frac{1}{k+2s-1}
\end{equation}
which are the ``shifted'' harmonic numbers. In particular, for $s=1/2$ one recovers the standard harmonic numbers,
which explains the name of the process.
}
\end{remark}

For a system of size $N$ and reservoirs parameters $0\le\rho_l\le \rho_r < \infty$
we denote by $\mu_{N,\rho_l,\rho_r}$ the  invariant measure of the process
$\{\eta(t)\,,\,t\ge 0\}$ of Definition \ref{process}, i.e. the ``non-equilibrium steady state''
of the boundary-driven harmonic process with parameter $s>0$.
To alleviate the notation we do not write in the measure the dependence
on the parameter $s$.

As a particular case, in the equilibrium set-up $\rho_l=\rho_r$,
one can check that the harmonic process with parameter $s>0$
has a reversible invariant measure given by a product of
Negative Binomial distributions with shape parameter $2s$ and mean $2s \theta$.
Namely, considering the univariate probability mass function
\be
\label{negbin-ew}
\nu_\theta(n):=\frac{1}{n!}\frac{\Gamma(2s+n)}{\Gamma(2s)}\left(\frac{\theta}{1+\theta}\right)^n \left(\frac1{1+\theta}\right)^{2s}
\qquad\qquad n\in\mathbb{N}_0 , \quad \theta \ge0
\ee
with mean
$$
\sum_{n=0}^\infty n \nu_\theta(n) = 2s \theta,
$$
and defining the product law
\be
\label{rev-mes}
\mu_{N,\rho_l,\rho_l}({\eta}) :=  \prod_{i=1}^N \nu_{\rho_l}(\eta_i) \qquad\qquad   \eta\in\Omega_N , \quad \rho_l >0
\ee
then one has
$\langle f, {\mathscr L}g \rangle = \langle {\mathscr L}f, g \rangle$,
where $\langle \cdot,	\cdot \rangle$ denotes the scalar product in the Hilbert space $L^2(\mathbb{N}^N, \mu_{N,\rho_l,\rho_l})$.

In the non-equilibrium case ($0\le\rho_l < \rho_r<\infty$) the stationary measure was computed in \cite{fg}
by a combination of stochastic duality and quantum inverse scattering method. Define the (scaled) factorial moment
of order $\xi =(\xi_1, \ldots, \xi_N)\in \mathbb N_0^N$ as
\beq\label{eq:Gasmu}
 G(\xi) =  \sum_{\eta \in \mathbb N_0^N} \mu_{N,\rho_l,\rho_r}(\eta) \left[ \prod_{i=1}^N \frac{\eta_i!}{(\eta_i-\xi_i)! }\cdot \frac{\Gamma(2s)}{\Gamma(2s+\xi_i)} \right] \;.
\eeq
Then the following result is available:
\bt[Factorial moments, \cite{fg}]
\label{SFM-new}
Using the notation $|\eta| = \sum_{i=1}^N \eta_i$,
the scaled factorial moments of the non-equilibrium steady state are given by
\be
\label{formula-new}
 G(\xi)  = \sum_{\eta \in \mathbb N_0^N}\rho_r^{|\xi|-|\eta|}(\rho_\ell-\rho_r)^{|\eta|} \, \prod_{i=1}^N \; \binom{\xi_i}{\eta_i}\; f_{i}(\eta)
\ee
with
\be
\label{fi-new}
f_{i}(\eta):= \prod_{j=1}^{\eta_i} \; \frac{2s(N+1-i)-j+\mathcal{N}^+_i(\eta)}{2s(N+1)-j+\mathcal{N}^+_i(\eta)}
\qquad\text{and} \qquad {\mathcal N}^+_i(\eta):= \sum_{k=i}^N \eta_k .
\ee
\et
The  steady state of the boundary driven harmonic process can the be reconstructed in terms of the factorial moments \eqref{formula-new} via
the formula
\begin{equation}\label{eq:measureee}
\mu(\eta) = \sum_{\xi \ge \eta} G(\xi)  \Big[\prod_{i=1}^N \frac{(-1)^{\xi_i-\eta_i}}{\xi_i!} {\xi_i \choose \eta_i} \frac{\Gamma(2s+\xi_i)}{\Gamma(2s)} \Big]\,.
\end{equation}

\medskip
\noindent

\section{The non-equilibrium steady state}
\label{pgf-section}

In this section we identify the non-equilibrium steady state of the harmonic model  in \eqref{eq:measureee}
as a mixture measure. In the equilibrium set-up ($\rho_l=\rho_r$) the invariant
measure is reversible and is an homogeneous (Gibbs) product measure.
In non-equilibrium  ($\rho_l\neq\rho_r$) we shall prove that
the invariant measure is a {\em mixture} of inhomogeneous product measures. The mixing measure is related to the order statistics of
uniform i.i.d. random variables when $2s$ is an integer,
and more generally to the ``ordered Dirichlet distribution''
when $2s$ is not an integer.

\subsection{Preliminaries: order statistics of uniform random variables}
In the following lemmata, we recall  a few facts about the order statistics of i.i.d. uniforms
on the unit interval. See \cite{abn}, \cite{order} for more details.

\bl[Marginals]
\label{lemma1-order}
Let $U_1, \ldots, U_n$ denote $n$ independent uniforms on $[0,1]$ and
denote their ascen\-ding order statistics by $U_{1,n}\leq U_{2,n}\leq \ldots \leq  U_{n,n}$ .
Let $1\le n_1 \le n$ then the marginal probability density of the $U_{n_1,n}$ is
\be
f_{U_{n_1,n}}(u_1) = \frac{n!}{(n_1-1)!(n-n_1)!} \:\cdot u_1^{n_1-1} (1-u_1)^{n-n_1} \cdot \indic{0\le u_1 \le 1}.
\ee
For a given $1\le k\le n$ this generalizes as follows:
if  $1 \le  n_1 <\ldots < n_k   \le  n$ then the joint probability density of
$(U_{n_1,n}, U_{n_2,n}, \ldots, U_{n_k,n})$ is
\beq\label{orderdirichletdens}
f_{(U_{n_1,n}, U_{n_2,n}, \ldots, U_{n_k,n})}(u_1, u_2,\ldots, u_k)
= n! \left[\prod_{i=1}^{k+1} \frac{(u_i-u_{i-1})^{n_i-n_{i-1}-1}}{(n_i-n_{i-1}-1)!} \right]
\,\indic{0\leq u_1 \leq \ldots \leq u_k \leq 1}
\eeq
where we used the convention $n_0=0$, $n_{k+1}=n+1$, $u_0=0$ and $u_{k+1}=1$.
\el
It is easy to see that the sequence of order statistics of continuous random
variables is Markov.
\bl[Markov property]
\label{lemma2-order}
Let $U_1, \ldots, U_n$ denote $n$ independent uniforms on $[0,1]$ and
denote their ascen\-ding order statistics by $U_{1,n}\leq U_{2,n}\leq \ldots \leq  U_{n,n}$ .
Then the order statistics forms a Markov chain, i.e. for all $1\le m \le n$, the sets of order statistics
$(U_{1,n},\ldots, U_{m-1,n})$ and $(U_{m+1,n},\ldots,U_{n,n})$
become conditionally independent if $U_{m,n}$ is fixed. Therefore
for the joint densities we may write
\begin{eqnarray}
&& f_{U_{1,n},\ldots, U_{m-1,n}, U_{m+1,n},\ldots, U_{n,n}\,|\, U_{m,n} }(u_1, \ldots,u_{m-1},u_{m+1},\ldots,u_n\,|\,u_m ) =\\ \nn
&& \hspace{1.cm} f_{U_{1,n},\ldots, U_{m-1,n}\,|\, U_{m,n}}(u_1, \ldots,u_{m-1}\,|\,u_m) \cdot
f_{{U}_{m+1,n},\ldots, {U}_{n,n}\,|\, U_{m,n}}(u_{m+1},\ldots,u_n\,|\,u_m).
\end{eqnarray}
\el
We also have the following important result: the conditional distribution of the order statistics (conditioned
on another order statistic) is related to the distribution of order statistics from a (smaller) population
whose distribution function is a truncated form of the original distribution function.

\bl[Left/right truncation]
\label{lemma3-order}
Let $U_1, \ldots, U_n$ denote $n$ independent uniforms on $[0,1]$ and
denote their ascen\-ding order statistics by $U_{1,n}\leq U_{2,n}\leq \ldots \leq  U_{n,n}$ .
Then, for $1\le m \le n$ and $u_m\in (0,1)$, the conditional distribution of $(U_{1,n},\ldots, U_{m-1,n})$, given that $U_{m,n} = u_m$
is the same as the distribution of the order statistics $(U^{\star}_{1,m-1},\ldots, U^{\star}_{m-1,m-1})$ obtained from a
sample of size $m-1$ from a population whose distribution is
uniform on $[0,u_m]$,
i.e.
\be
 f_{U_{1,n},\ldots, U_{m-1,n}\,|\, U_{m,n}}(u_1, \ldots,u_{m-1}\,|\,u_m ) = f_{U^{\star}_{1,m-1},\ldots, U^{\star}_{m-1,m-1}}(u_1, \ldots,u_{m-1}) \;.
\ee
Similarly, the conditional distribution of $(U_{m+1,n},\ldots,U_{n,n})$, given that $U_{m,n} = u_m$
is the same as the distribution of the order statistic $(\tilde{U}_{1,n-m},\ldots,\tilde{U}_{n-m,n-m})$ obtained from a
sample of size $n-m$ from a population whose distribution is
uniform on $[u_m,1]$,
i.e.
\be
f_{U_{m+1,n}\ldots, U_{n,n}\,|\, U_{m,n} }(u_{m+1},\ldots,u_n\,|\,u_m ) =
f_{\tilde{U}_{1,n-m},\ldots, \tilde{U}_{n-m,n-m}}(u_{m+1},\ldots,u_n) \;.
\ee
\el
Combining together Lemma \ref{lemma2-order} and Lemma \ref{lemma3-order} we obtain
the following property for the conditional distribution of the order statistics of i.i.d. uniform random variables
on the interval $[0,1]$.
\bl[Conditional distribution]\label{lemma4-order}
With the same hypotheses and notations of Lemma \ref{lemma2-order} and Lemma \ref{lemma3-order}
we have
\begin{eqnarray}
\label{main-order}
&& f_{U_{1,n},\ldots, U_{m-1,n}, U_{m+1,n}\ldots, U_{n,n}\,|\, U_{m,n} }(u_1, \ldots,u_{m-1},u_{m+1},\ldots,u_n\,|\,u_m ) =\\ \nn
&& \hspace{1.cm} f_{U^{\star}_{1,m-1},\ldots, U^{\star}_{m-1,m-1}}(u_1, \ldots,u_{m-1}) \cdot
f_{\tilde{U}_{1,n-m},\ldots, \tilde{U}_{n-m,n-m}}(u_{m+1},\ldots,u_n) \;.
\end{eqnarray}
\el
Finally, we will use the following large deviation result for the
sample paths of the order statistics.
Let $U_1,\ldots, U_n$ be a random i.i.d. sample from a uniform distribution on $[0,1]$,
and let $U_{1,n}\leq U_{2,n}\leq \ldots \leq  U_{n,n}$ denote the order statistics obtained from this sample.
Using the convention $U_{n+1,n}:=1$, we define the sample path of the order statistics by
$$
U_n(t)= U_{\lfloor (n+1)t\rfloor +1,n} \qquad\qquad \text{for all } t\in[0,1]
$$
where $\lfloor y \rfloor$ denotes the largest integer that is smaller or equal to $y$.
Then we have the following functional  {Large Deviation Principle (LDP)} for the sample paths of the order statistics.
\bl[Sample path large deviation for order statistics, \cite{duffy-macci-torrisi}]
\label{ld-order}
Let $D[0,1]$ denote the space of c\`{a}dl\`{a}g functions
on the unit interval, equipped with Skorohod topology.
Let $A_{0,1} \subset D[0,1]$ denote the closed set of non-decreasing functions
$f:[0,1]\to \mathbb{R}$ such that $f(x)\ge 0$ and $f(1)=1$.
Then the sample paths $U_n(\cdot)$ satisfy the large deviation principle
with rate function
$$
 J(u) =
\left\{\begin{array}{ll}
-\int_{0}^1 \log(u'(x)) dx  \qquad&  \text{if } u \in A_{0,1} \text{ is strictly increasing}\\
\infty & \text{otherwise}
 \end{array}\right.
$$
\el

\subsection{Stationary measure as a probabilistic mixture.}\label{Mixt}
\bt[Mixture structure of the NESS]
\label{mixtmeasthm}
Let $2s\in \mathbb{N}$ and $N\in \mathbb{N}$ and assume without
loss of generality that $0\le \rho_l \le  \rho_r < \infty$. Define $$n:= 2s(N+1)-1$$ and let $U_1,\ldots, U_n$ be independent random variables
with common uniform distribution on the interval $(0,1)$.
Consider  the  distribution of the $N$-dimensional vector $(U_{2s,n}, U_{4s,n}, \ldots, U_{2sN,n})$
obtained as a marginal of the order statistics $U_{1,n} \le  U_{2,n} \le \ldots \le  U_{n,n}$,
whose probability density reads (using \eqref{orderdirichletdens} with $k=N$
and $n_i = 2s i$ for $i=1,\ldots,N$)
\beq\label{2sdirdens}
f_{(U_{2s,n}, U_{4s,n}, \ldots, U_{2sN,n})}(u_1, \ldots, u_N)
= \frac{\Gamma(2s(N+1))}{\Gamma (2s)^{N+1}}\cdot \prod_{i=1}^{N+1} (u_i-u_{i-1})^{2s-1} \cdot
\,\indic{0 \leq u_1 \leq \ldots \leq u_N \leq 1}
\eeq
with the convention $u_0=0$ and $u_{N+1}=1$.
Then, the non-equilibrium steady state of the open harmonic process
of Definition \ref{process}
is equal to
\be\label{mixtresult}
\mu_{N,\rho_l, \rho_r}(\eta)= \mathbb{E}\left( \prod_{i=1}^N \nu_{\Theta_{2si,n}} (\eta_i)\right)
\ee
where $\nu_\theta$ is the Negative Binomial law defined in \eqref{negbin-ew} and
the expectation $\mathbb{E}$ is w.r.t. the random variables $(\Theta_{2s,n}\ldots,\Theta_{2sN,n})$
obtained as a marginal of the order statistics $\Theta_{1,n} \le   \ldots \le  \Theta_{n,n}$
of the independent random variables
\be
\label{rho-i}
\Theta_{i}= \rho_l + (\rho_r-\rho_l) U_{i} \qquad\qquad i = 1,\ldots,n 
\ee
that have  uniform  distribution on $[\rho_l,\rho_r]$.
More explicitely
\begin{eqnarray}\label{mixtresult-explicit}
\mu_{N,\rho_l, \rho_r}(\eta)
&= &
\frac{\Gamma(2s(N+1))}{\Gamma(2s)^{N+1}}\frac{1}{(\rho_r-\rho_l)^{2s(N+1)-1}}\cdot
 \int_{\rho_l}^{\rho_r}d\theta_1 \int_{\theta_1}^{\rho_r}d\theta_2\cdots  \int_{\theta_{N-1}}^{\rho_r}d\theta_N  \prod_{i=1}^{N+1} (\theta_{i}-\theta_{i-1})^{2s-1}   \nonumber  \\
& & 
\prod_{i=1}^{N}\frac{1}{\eta_i!}\frac{\Gamma(2s+\eta_i)}{\Gamma(2s)}\left(\frac{ \theta_i}{1+ \theta_i}\right)^{\eta_i} \left(\frac1{1+ \theta_i}\right)^{2s}
\end{eqnarray}
with the convention $\theta_0=\rho_l$ and $\theta_{N+1}=\rho_r$.
\et

\medskip
\begin{remark}[The case of non-integer $2s$]
{\em
The integral representation of the stationary measure given  in \eqref{mixtresult-explicit} has a meaning
even when one drops the assumption of $2s$ being an integer.
In this case, the law of the mixing measure is related to the ``ordered Dirichlet distribution''.
More precisely, for general $s>0$ we have
\be\label{mixtresult2}
\mu_{N,\rho_l, \rho_r}(\eta)= \mathbb{E}\left( \prod_{i=1}^N \nu_{S_i} (\eta_i)\right)
\ee
where now the expectation
$\mathbb{E}$ denotes expectation w.r.t. the joint distribution of the
random variables $(S_1,\ldots,S_N)$ defined by
\begin{equation*}
S_i= \rho_l+ (\rho_r-\rho_l) V_{i}
\end{equation*}
where $(V_1,\ldots,V_N)$  is the random vector
with joint probability density
\begin{equation*}
f_{(V_1,\ldots,V_N)}(v_1, \ldots, v_N)
=\frac{\Gamma(2s(N+1))}{\Gamma (2s)^{N+1}} \cdot  \prod_{i=1}^{N+1} (v_{i}-v_{i-1})^{2s-1} \cdot \indic{v_0=0\leq v_1\leq v_2\leq\ldots v_N\leq v_{N+1}=1}\;.
\end{equation*}
Such distribution arises from the sum of the components of the symmetric Dirichlet distribution.
Indeed, let $R_i:= V_i-V_{i-1}$ for $ i=1,\ldots, N+1$ with $V_0 = 0$ and $V_{N+1}=1$, then its inverse transformation is
\begin{equation*}
V_i= \sum_{j=1}^i R_j \qquad \text{ for } i=1,\ldots, N.
\end{equation*}
The joint distribution of $(R_1, \ldots, R_{N+1})$ reads
\begin{equation*}
f_{(R_1,\ldots,R_{N+1})}(r_1, \ldots, r_{N+1})
= \frac{\Gamma(2s(N+1))}{\Gamma (2s)^{N+1}} \prod_{i=1}^{N+1} r_i^{2s-1}\indic{\Sigma_{N+1}}(r_1,\ldots,r_{N+1})\nn
\end{equation*}
which is the Dirichlet distribution on the $(N+1)$-dimensional simplex
\[
\Sigma_{N+1}=\{(r_1, \ldots, r_{N+1}): 0\leq r_i\leq 1\ \text{for all}\ i,\, r_1+\ldots +r_{N+1}=1\}
\]
with all parameters equal to $2s>0$.
In the case of integer $2s$, one recovers the representation \eqref{mixtresult}
from the representation \eqref{mixtresult2} using
$$
(R_1,R_2,\ldots,R_{N+1}) =(U_{2s,n}-U_{0,n}, U_{4s,n}-U_{2s,n}, \ldots, U_{2s(N+1),n}- U_{2sN,n})
$$
with  the convention $U_{0,n}=0$ and $U_{2s(N+1),n}=1$, which is
the well-known relation between the symmetric Dirichlet distribution
with parameter $2s$ on the $(N+1)$-dimensional simplex and the vector
constructed from differences (with gaps $2s$) of the order statistics of $n=2s(N+1)-1$
 i.i.d. uniform random variables on the unit interval.
}
\end{remark}

\subsection{Proof of Theorem \ref{mixtmeasthm}}
In this section we provide a proof of Theorem~\ref{mixtmeasthm}.
We also refer the reader to Appendix~\ref{app:compare}  where it is shown that the integral representation \eqref{mixtresult-explicit} is identical to the closed-form expression in \eqref{eq:measureee}.

\subsubsection{Moment generating function}
The strategy to prove Theorem \ref{mixtmeasthm} is to use the moment generating function
to characterize the stationary measure.
Define the set
$$
{{\mathcal A}_{N,\rho_l,\rho_r} = \left\{\mathsf{h}=(\mathsf{h}_1,\ldots,\mathsf{h}_N)\in \mathbb R^N \,:\, |\mathsf{h}_i| \le \log \left( 1 + \frac{1}{\rho_r} \right)  \quad \text{for }i=1,\ldots,N\right\}}.
$$
For  $\mathsf{h} \in {\mathcal A}_{N,\rho_l,\rho_r} $, let us denote by $\Psi_{N,\rho_l,\rho_r}(\mathsf{h})$ the moment generating function (MGF) of the
non-equilibrium steady state, i.e.
\be
\label{pgf-N}
\Psi_{N,\rho_l,\rho_r}(\mathsf{h})= \sum_{\eta}  \mu_{N,\rho_l,\rho_r}(\eta)  \prod_{i=1}^N  e^{\mathsf{h}_i \eta_i}.
\ee
Starting from the factorial moments \eqref{formula-new} we will compute the generating function
and show it coincides with the one of the law \eqref{mixtresult}.
We split the computation of the moment generating function into three steps,
which are given in Proposition \ref{pgf-unnested}, Proposition \ref{prop-unnestedintegrals}
and in Proposition \ref{prop-conclude}.

\subsubsection{N-fold sums}\label{Nfoldsumssubsection}
In this section  {we show} that the moment generating function $\Psi_{N,\rho_l,\rho_r}(\mathsf{h})$ can be written,  modulo multiplication by a factor,  as the composition of a function $\Phi_N: \mathbb R^N \to \mathbb R$  and the map
\begin{align} \label{ci-fct}
c_{N,\rho_\ell,\rho_r} : \mathbb R^N & \longrightarrow \mathbb R^N \\
(\mathsf{h}_1, \ldots, \mathsf{h}_N) & \longrightarrow \left( \frac {(\rho_r-\rho_\ell)\left(1-e^{\mathsf{h}_1} \right)}{1+\rho_r(1-e^{\mathsf{h}_1})}, \ldots, \frac {(\rho_r-\rho_\ell)\left(1-e^{\mathsf{h}_N} \right)}{1+\rho_r(1-e^{\mathsf{h}_N})} \right) \nn
\end{align}
i.e. the $i$-th component of the vector $c_{N,\rho_\ell,\rho_r}(\mathsf{h})$ is given by
\be\label{ci-fcti}
(c_{N,\rho_\ell,\rho_r}(\mathsf{h}))_i = c_{\rho_r, \rho_\ell,i}(\mathsf{h}_i) :=  \frac {(\rho_r-\rho_\ell)\left(1-e^{\mathsf{h}_i} \right)}{1+\rho_r(1-e^{\mathsf{h}_i})}.
\ee
We will see that the function $\Phi_N$ for which we will obtain an explicit formula in terms of an $N$-fold sum, does not depend on the boundary densities $\rho_l$ and $\rho_r$. The dependence on this parameters is then completely offloaded onto the map   $c_{N,\rho_\ell,\rho_r}$.
\bp[MGF, un-nested sums ]
\label{pgf-unnested}
For $\mathsf{h} \in \mathcal{A}_{N,\rho_l,\rho_r}$ we have that
\beq
\label{Phi-fct}
  \Psi_{N,\rho_l,\rho_r}(\mathsf{h})   = \, \prod_{i=1}^N \; \big(1+\rho_r(1- e^{\mathsf{h}_i})\big)^{-2s} \cdot \Phi_{N}(c_{N,\rho_\ell,\rho_r}(\mathsf{h}))
\eeq
with $c_{N,\rho_\ell,\rho_r}:\mathbb R^N \to \mathbb R^N$ defined in \eqref{ci-fct}-\eqref{ci-fcti} and
\begin{equation}
\label{Phi-unnested-sums}
 \Phi_{N}(c)=  \frac{\Gamma(2s(N+1))}{\Gamma(2s)} \sum_{\eta \in \mathbb N_0^N} \prod_{i=1}^N
 c_i^{\eta_i}\cdot
\frac{1}{\eta_i!}\frac{\Gamma(\eta_i+2s)} {\Gamma(2s)} \cdot \;
 \frac{\Gamma(2s(N+1-i)+\mathcal{N}^+_i(\eta))}{\Gamma(2s(N+2-i)+\mathcal{N}^+_{i}(\eta))}.
 \end{equation}
\ep
\bpr
The moment generating function can be rewritten in terms of the
scaled factorial moments  as follows:
\beq
\Psi_{N,\rho_l,\rho_r}(\mathsf{h})
&=& \sum_{\eta}\left[ \prod_{i=1}^N  \sum_{\xi_i=0}^{\eta_i} \binom{\eta_i}{\xi_i}(e^{\mathsf{h}_i}-1)^{\xi_i}\right] \mu_N(\eta)\nn \\
&=&\sum_{\xi} \left[ \prod_{i=1}^N \frac{1}{\xi_i!} \frac{\Gamma(2s+\xi_i)}{\Gamma(2s)}   (e^{\mathsf{h}_i} -1)^{\xi_i}\right]  G(\xi) \nn
\eeq
where it has been used that ${\eta_i \choose \xi_i} = 0$ for natural numbers $\xi_i > \eta_i$.
Therefore, as a consequence of Theorem \ref{SFM-new} we have
\beq
\Psi_{N,\rho_l,\rho_r}(\mathsf{h})
&=&
\sum_{\xi} \left[ \prod_{i=1}^N  \frac{\Gamma(2s+\xi_i)}{\Gamma(2s)\cdot \xi_i!} \,  \left(e^{\mathsf{h}_i}-1 \right)^{\xi_i} \right]
\; \sum_{\eta}\rho_r^{|\xi|-|\eta|}(\rho_\ell-\rho_r)^{|\eta|} \, \prod_{i=1}^N \; \binom{\xi_i}{\eta_i}\; f_{i}(\eta)
  \nn\\
&=&
\sum_{\substack{{\xi,\eta}\\{\eta\le \xi}}}   \prod_{i=1}^N \
\rho_r^{\xi_i-\eta_i}(\rho_\ell-\rho_r)^{\eta_i}
\binom{\xi_i}{\eta_i} \, \,  \; f_{i}(\eta) \left(e^{\mathsf{h}_i} - 1 \right)^{\xi_i} \frac{\Gamma(2s+\xi_i)}{\Gamma(2s)\cdot \xi_i!}
 \nn
\eeq
where we  used the notation $\eta\le \xi$ to indicate  that $\eta_i\le \xi_i$ for all $i\in\{1, \ldots, N\}$.
By exchanging the order of summations we obtain
\beq
\Psi_{N,\rho_l,\rho_r}(\mathsf{h})&=& \sum_{\eta} \prod_{i=1}^N (\rho_\ell-\rho_r)^{\eta_i}  \left(e^{\mathsf{h}_i}-1 \right)^{\eta_i} \; f_{i}(\eta) \sum_{\xi_i\ge \eta_i}  \binom{\xi_i}{\eta_i} \,     \rho_r^{\xi_i-\eta_i} \left(e^{\mathsf{h}_i}-1 \right)^{\xi_i-\eta_i} \,   \frac{\Gamma(2s+\xi_i)}{\Gamma(2s)\cdot \xi_i!}.\nn
\eeq
The sum of the $\xi$ variables can now be performed
using that for all $i\in\{1,\ldots,N\}$
\beq
\sum_{\xi_i\ge \eta_i}  \binom{\xi_i}{\eta_i} \,     \rho_r^{\xi_i-\eta_i} \left(e^{\mathsf{h}_i}-1 \right)^{\xi_i-\eta_i} \frac{\Gamma(2s+\xi_i)}{\Gamma(2s)\cdot \xi_i!}
 &=& \frac{\Gamma(\eta_i+2s)} {\Gamma(2s) \cdot \eta_i!}  \sum_{\xi_i\ge \eta_i}   \frac{\Gamma(2s+\xi_i)}{\Gamma(\eta_i+2s)\cdot (\xi_i-\eta_i)!} \,   \rho_r^{\xi_i-\eta_i} \left(e^{\mathsf{h}_i}-1 \right)^{\xi_i-\eta_i} \nn\\
 &=&
 {\color{black} \frac{\Gamma(\eta_i+2s)} {\Gamma(2s) \cdot \eta_i!} \sum_{k_i \ge 0}  \frac{\Gamma(2s+\eta_i+k_i)}{\Gamma(2s+\eta_i) \cdot k_i!}    \,     (\rho_r(e^{\mathsf{h}_i}-1))^{k_i}}=\nn\\
  &=&   \frac{\Gamma(\eta_i+2s)} {\Gamma(2s) \cdot \eta_i!}\,  \frac 1 {(1-\rho_r(e^{\mathsf{h}_i}-1))^{\eta_i+2s}}.\nn
\eeq
where in the last equality we have used the identity
\be
\label{identity-bin-sum}
 \frac{1}{(1-x)^{a}}=\sum_{k=0}^\infty \frac{\Gamma(a+k)}{\Gamma(a) \cdot k!} x^k \qquad\qquad |x| < 1.
\ee
Thus we arrive to
\beq
\Psi_{N,\rho_l,\rho_r}(\mathsf{h})&=& \sum_{\eta} \prod_{i=1}^N (\rho_\ell-\rho_r)^{\eta_i}  \left(e^{\mathsf{h}_i} - 1 \right)^{\eta_i} f_{i}(\eta)\;
 \frac{\Gamma(\eta_i+2s)} {\Gamma(2s) \eta_i!}\,     \frac {1} {(1-\rho_r(e^{\mathsf{h}_i}-1))^{\eta_i+2s}}.\nn
 \eeq
Equivalently, multiplying both sides by $\prod_{i=1}^N \; \big(1+\rho_r(1-e^{\mathsf{h}_i})\big)^{2s}$
we rewrite this identity in terms of the function $\Phi_N$  defined in \eqref{Phi-fct} as
\beq
\label{capital}
\Phi_{N}(c)
&=& \sum_{\eta}\prod_{i=1}^N  c_i^{\eta_i}\,\frac{\Gamma(\eta_i+2s)} {\Gamma(2s) \cdot \eta_i!} \cdot f_{i}(\eta)
\eeq
with $c_i$ as given in \eqref{ci-fct}. Recalling the definition of the functions $f_{i}$ in \eqref{fi-new} and
using the convention ${\mathcal N}_{N+1}^+(\eta)=0$, we write $\prod_{i=1}^N f_{i}(\eta)$
as a telescopic product
\be
\nn
\prod_{i=1}^N f_{i}(\eta)= \prod_{i=1}^N \prod_{k={\mathcal N}_{i+1}^+(\eta)}^{\mathcal{N}_{i}^+(\eta)-1} \; \frac{2s(N+1-i)+k}{2s(N+1)+k} .
\ee
As a consequence
\beq
\prod_{i=1}^N f_{i}(\eta)&=&\frac {\Gamma(2s(N+1))}{\Gamma(2s(N+1)+\mathcal{N}^+_1(\eta))} \cdot \prod_{i=1}^N
\frac {\Gamma(2s(N+1-i)+\mathcal{N}^+_i(\eta))}{\Gamma(2s(N+1-i)+\mathcal{N}^+_{i+1}(\eta))}  \nn\\
&=& \frac {\Gamma(2s(N+1))}{\Gamma(2s)} \cdot \prod_{i=1}^N \frac{\Gamma(2s(N+1-i)+\mathcal{N}^+_i(\eta))}{\Gamma(2s(N+1-(i-1))+\mathcal{N}^+_i(\eta))}.\nn
\eeq
Inserting this last expression in \eqref{capital}, the result of the proposition follows.
\epr
\br[MGF, nested sums]
\label{remark-nested-discrete}
{\em
There is a one-to-one relation between the set of configurations $\eta\in \mathbb N_0^N$ and the
set of $N$-tuples $\{(m_1, \ldots, m_N)\in \mathbb N_0^{N}: m_1\ge m_2  \ge \ldots \ge m_N\ge 0\}$.
This implies that the moment generating function can also be written
as nested sums.
Then we have
\beq
\Phi_{N}(c)&=&  \frac{\Gamma(2s(N+1))}{\Gamma(2s)}   \sum_{m_1 \ge \ldots \ge m_N\ge 0}\prod_{i=1}^N  c_i^{m_i-m_{i+1}}\, \frac{\Gamma(m_i-m_{i+1}+2s)} {\Gamma(2s) (m_i-m_{i+1})!} \cdot
\frac{\Gamma(2s(N+1-i)+m_i)}{\Gamma(2s(N+2-i)+m_i)}\nn
\eeq
with the convention $m_{N+1}=0$.
This easily follows from Proposition \ref{pgf-unnested} by implementing the change of variables:
$$
\eta=(\eta_1, \ldots, \eta_N) \longrightarrow m =(m_1, \ldots, m_N), \qquad \text{with} \qquad m_i:= \mathcal{N}_i^+(\eta)
$$
from which one has
$
\eta_i(m)=m_i-m_{i+1}.
$
}
\er

\subsubsection{N-fold integrals} \label{N-fold integrals}
We proceed further by moving from a representation of the moment generating function
with $N$ sums to one involving $N$ integrals. This will be useful to recognize the invariant
distribution of the harmonic process as a mixture.
\bp[MGF, un-nested integrals ]
\label{prop-unnestedintegrals}
We have
\begin{equation}
\label{Phi-unnested-integrals}
\begin{split}
 \Phi_{N}(c)
 &=\frac{\Gamma(2s(N+1))}{\Gamma(2s)^{N+1}}
 \int_{0}^{1}dt_1 \cdots  \int_{0}^{1}dt_N\prod_{i=1}^N
  t_i^{2s(N-i+1)-1} (1-t_i)^{2s -1} \left(\frac{1}{1-c_i \prod_{j=1}^i t_j}\right)^{2s}.
 \end{split}
\end{equation}
\ep
\bpr
We prove that \eqref{Phi-unnested-integrals} coincides with \eqref{Phi-unnested-sums}
using again the identity \eqref{identity-bin-sum}.
Indeed, plugging this identity in \eqref{Phi-unnested-integrals} we have
\begin{equation}
\begin{split}
 \Phi_{N}(c)
 &=\frac{\Gamma(2s(N+1))}{\Gamma(2s)^{N+1}}
 \int_{0}^{1}dt_1 \cdots  \int_{0}^{1}dt_N\prod_{i=1}^N  t_i^{2s(N-i+1)-1} (1-t_i)^{2s -1}
 \sum_{\eta_i=0}^{\infty} \frac{\Gamma(2s+\eta_i)}{\Gamma(2s)\eta_i!}
 \left(c_i \prod_{j=1}^i t_j\right)^{\eta_i}.\nn
 \end{split}
\end{equation}
Collecting the powers of $t_i$ and recalling the definition $\mathcal{N}^+_i(\eta) = \sum_{k= i}^N \eta_k$
this can be rewritten as
\begin{equation}
 \Phi_{N}(c)
 =\frac{\Gamma(2s(N+1))}{\Gamma(2s)^{N+1}}
 \sum_{\eta}  \prod_{i=1}^N   \frac{\Gamma(2s+\eta_i)}{\Gamma(2s)\eta_i!}
\cdot  c_i^{\eta_i}    \int_{0}^{1}
  t_i^{2s(N-i+1)+\mathcal{N}_i^+(\eta)-1} (1-t_i)^{2s -1} dt_i.
\nn
\end{equation}
Using that for all $a,b >0$
$$
\int_{0}^1  x^{a-1} (1-x)^{b-1} dx = \frac{\Gamma(a)\Gamma(b)}{\Gamma(a+b)}
$$
it then follows
\begin{equation}
 \Phi_{N}(c)
 =\frac{\Gamma(2s(N+1))}{\Gamma(2s)^{N+1}}
 \sum_{\eta}  \prod_{i=1}^N   c_i^{\eta_i} \: \frac{\Gamma(2s+\eta_i)}{\Gamma(2s)\eta_i!}
\cdot
 \frac{\Gamma(2s(N+1-i)+\mathcal{N}_i^+(\eta))\cdot \Gamma(2s)}{\Gamma(2s(N+2-i)+\mathcal{N}_i^+(\eta))}\nn
\end{equation}
which reproduces \eqref{Phi-unnested-sums} after simplifications.
\epr
\br[MGF, nested integrals ]
\label{prop-nested-integrals}
{\em
Similarly to the discrete case (see Remark \ref{remark-nested-discrete}),
one can also write an expression in terms of nested integrals.
We have
\beq
\label{Phi-nested-integrals}
\Phi_{N}(c)&=&
\frac{\Gamma(2s(N+1))}{\Gamma(2s)^{N+1}}\cdot
 \int_0^{1}du_1 \int_{u_1}^1du_2\cdots  \int_{u_{N-1}}^1du_N
\,\prod_{i=1}^{N+1} (u_{i}-u_{i-1})^{2s-1} \frac{1}{\big(1-c_i(1-u_i)\big)^{2s}} \nn\\
\eeq
where we recall the convention $u_{0}=0$ and $u_{N+1}=1$.
The result easily follows from Proposition \ref{prop-unnestedintegrals} by implementing the change
of variables
$ u_i = 1-\prod_{j=1}^i t_j.
$
Inverting this mapping one gets
$$
t_i = \frac{1-u_{i}}{1-u_{i-1}} \qquad\qquad\text{and}\qquad\qquad 1-t_i = \frac{u_{i}-u_{i-1}}{1-u_{i-1}}
$$
which substituted in \eqref{Phi-unnested-integrals} yields \eqref{Phi-nested-integrals}.
}
\er

\subsubsection{Concluding the proof}

The last step in the proof of Theorem  \ref{mixtmeasthm} consists in recognizing
in the expression \eqref{Phi-nested-integrals} the probability
generating function of the probability measure \eqref{mixtresult}.
We recall that the moment generating function of
a Negative Binomial distribution with law \eqref{negbin-ew} is given by
\be
\label{pgf-negbin}
M_{\theta}(\mathsf{h}) = \sum_{n=0}^\infty e^{\mathsf{h}n} \nu_\theta(n)=\left(\frac1{1+\theta(1- e^\mathsf{h})}\right)^{2s} \qquad\qquad \text{ for } |\mathsf{h}|< \log \left( 1+ \tfrac1\theta \right) .
\ee

\bp[MGF, mixture]
\label{prop-conclude}
For $\mathsf{h} \in\mathcal{A}_{N,\rho_l,\rho_r}$ we have
\be
\label{mix-2s-int}
\Psi_{N,\rho_l,\rho_r}(\mathsf{h})
 = \mathbb{E}\left[\prod_{i=1}^N M_{\Theta_{2si,n}}(\mathsf{h}_i)\right]
\ee
where the expectation is w.r.t. the marginal distribution of the ascending order statistics
of the i.i.d. uniform random variables defined in  \eqref{rho-i}.
\ep
\bpr
We observe that using \eqref{ci-fct}, namely
$$
 c_i=\frac {(\rho_r-\rho_\ell)\left(1-e^{\mathsf{h}_i} \right)}{1+\rho_r(1-e^{\mathsf{h}_i})}
$$
we have
$$
\frac{1}{1-c_i (1-u_i)}= \frac{\big(1+\rho_r(1-e^{\mathsf{h}_i})\big)}{1+ \big(\rho_l + (\rho_r-\rho_l) u_i \big)(1-e^{\mathsf{h}_i})}
$$
Inserting this into \eqref{Phi-nested-integrals}
and recalling the relation \eqref{Phi-fct},
the moment generating function of the
non-equilibrium steady state is given by
\beq
\label{genfunct-noneq}
\begin{split}
\Psi_{N,\rho_l,\rho_r}(\mathsf{h})&=
\frac{\Gamma(2s(N+1))}{\Gamma(2s)^{N+1}}
 \int_0^{1}du_1 \int_{u_1}^1du_2\cdots  \int_{u_{N-1}}^1du_N  \\
 &
\hspace{4.cm}\prod_{i=1}^{N+1} (u_{i}-u_{i-1})^{2s-1} \cdot \prod_{i=1}^N\frac{1}{\big(1+ \big(\rho_l + (\rho_r-\rho_l) u_i \big)(1-e^{\mathsf{h}_i})\big)^{2s}} .
\end{split}
\eeq
Therefore, using  \eqref{2sdirdens}, \eqref{rho-i} and \eqref{pgf-negbin} we obtain \eqref{mix-2s-int}.
\epr

\section{Pressure}

In this section we use the  {characterization} of the stationary measure  {in Theorem \ref{Mixt}} to compute the pressure
associated to the non-equilibrium steady state. We will reproduce the expression predicted by
the Macroscopic Fluctuation Theory by first conditioning to a given realization of the random
local parameters and then using the large deviation properties of those local parameters.

\bt[Pressure]
\label{theo-pr}
Let $h:[0,1]\to \mathbb{R}$ be a smooth function. Define the pressure of the open symmetric harmonic process
 as
\be\label{PP}
{P}(h) :=  \lim_{N\to\infty} \frac{1}{N} \log \mathbb{E} \Big[e^{\sum_{i=1}^N \eta_i h(\frac{i}{N})} \Big] \,.
\ee
Then the pressure admits the following variational expression:
\be
{P} (h) =
 \sup_{\underset{\theta(1)=\rho_r}{\underset{\theta(0) = \rho_l}{\underset{\text{strictly increasing }}{\theta:[0,1]\to \mathbb{R}_+}}}}
 \Big[P(h,\theta) - J(\theta)\Big]
\ee
where
\be
\label{pressure-h-theta}
P(h,\theta)=  2s \int_{0}^1
 \log \Big(
 \frac{1}{1+(1-e^{h(x)})\theta(x)} \Big)
dx
\ee
and
\be
J(\theta) =
 - 2s \int_{0}^1 \log\Big(
 \frac{\theta'(x)}{\rho_r -\rho_l}\Big) dx.
\ee
\et

\bpr
Recalling  Proposition \ref{prop-conclude}, we have
\begin{eqnarray}\label{EEE}
\mathbb{E} \Big[e^{ \sum_{i=1}^N \eta_i h(\frac{i}{N})} \Big]
& = &
 \mathbb{E}\left[\prod_{i=1}^N M_{\Theta_{2si,n}}\left( h \left( \tfrac{i}{N}\right) \right) \right]\\
& = &
\mathbb{E}\left[\prod_{i=1}^N\left(\frac1{1+\Theta_{2si,n} \Big( 1-e^{h(\tfrac{i}{N})}\Big) }\right)^{2s}\right]
\end{eqnarray}
where $n=2s(N+1)-1$.
Introducing the sample path of the order statistics
$$
\Theta_{n}(x) = \Theta_{\lfloor (n+1)x\rfloor +1,n} \qquad x\in[0,1]
$$
with the convention $\Theta_{n+1,n}:=\rho_r$, we arrive to
\begin{eqnarray*}
\mathbb{E} \Big[e^{ \sum_{i=1}^N \eta_i h(\frac{i}{N})} \Big]
& = &
\mathbb{E}\left[\prod_{i=1}^N\left(\frac1{1+\Theta_{n}(\frac{2si}{n}) \Big( 1-e^{h(\frac{i}{N})} \Big)}\right)^{2s}\right]\\
& = &
\mathbb{E}\left[\exp\left\{2s\sum_{i=1}^N \left(\log\left(\frac1{1+\Theta_{n}(\frac{i}{N}) \Big( 1-e^{h(\frac{i}{N})} \Big) }\right)+o(1)\right) \right\}\right]
\end{eqnarray*}
where $o(1)$ to $0$ as $N\to\infty$, uniformly.

For an increasing function $\theta:[0,1]\to\mathbb{R}$ we define
$$
P_N(h, \theta) = \frac{2s}{N} \sum_{i=1}^N\log \left(\frac1{1+\theta(\frac{i}{N}) \Big( 1-e^{h(\frac{i}{N})} \Big) }\right) .
$$
By using the properties of conditional expectation, this allows to rewrite the generating function of the empirical distribution
as the conditional expectation of an exponential functional
\begin{eqnarray}
\label{exp-functional}
\mathbb{E} \Big[e^{ \sum_{i=1}^N \eta_i h(\frac{i}{N})} \Big]
& = &
\mathbb{E}\Big[\mathbb{E}\Big[\exp\left\{N(P_N(h, \Theta)+o(1))\right\} | \Theta\Big]\Big]
\end{eqnarray}
where we denote by $\Theta$ the collection of random variables $(\Theta_{n}(\frac{i}{N}))_{i=1,\ldots,N}$.
Observe that by Riemann approximation
$$
\lim_{N\to\infty}P_N(h,\theta) =  P(h,\theta) = 2s \int_{0}^1 \log \left(\frac1{1+\theta(x) \big( 1-e^{h(x)} \big) }\right) dx
$$
and recalling (see Lemma \ref{ld-order}) that the sample path of the order statistics satisfy the LDP with good rate function
$$
J(\theta) =
\left\{\begin{array}{ll}
-2s \int_{0}^1 \log(\frac{\theta'(x)}{\rho_r-\rho_l}) \, dx  \qquad&  \text{if } \theta \in A_{\rho_l,\rho_r} \text{ is strictly increasing}\\
\infty & \text{otherwise}
 \end{array}\right.
$$
the claim of the theorem follows by applying Varadhan's lemma to the exponentially growing functional
\eqref{exp-functional}.
\epr

 \section{Large deviations}

 In this section we prove that the sequence of empirical density measures $(L_N)_{N\ge 1}$
 satisfies a LDP.  One might think that knowing the pressure one could extract from it the
 large deviation function by using G\"artner-Ellis theorem. As we shall see and comment
 below this is not possible because the large deviation function is not convex.
 However we can obtain the large deviation function by following a direct approach
 that starts from the explicit knowledge of the (microscopic) stationary measure of the
 open harmonic model and proceed via a contraction principle.

\bt[Density large deviation]
\label{theo-ld}
The empirical profiles  of the open symmetric harmonic process
$$
L_N = \frac{1}{N} \sum_{i=1}^N \eta_i \delta_{\frac{i}{N}}
$$
satisfy a large deviation principle with good rate function
\be
\label{ldf}
 {I}(\rho) = \inf_{\underset{\theta(1)=\rho_r}{\underset{\theta(0) = \rho_l}{\underset{\text{strictly  increasing}}{\theta:[0,1]\to \mathbb{R}_+}}}}
 \Big[I(\rho,\theta) + J(\theta)\Big]
\ee
where
\be
\label{i-rho-theta}
I(\rho,\theta) = 2s \int_{0}^1 \left[
\frac{ \rho(x)}{2s} \log \frac{\rho(x)}{2s \theta(x)}
 + \left(1+\frac{\rho(x)}{2s}\right) \log\Big(\frac{1+ \theta(x)}{1+\frac{\rho(x)}{2s}}\Big)
 \right] dx
\ee
and
\be\label{Jt}
J(\theta) =
 - 2s \int_0^1 \log\Big(
 \frac{\theta'(x)}{\rho_r -\rho_l}\Big) dx \;.
\ee
\et
\noindent
Before proving the theorem we add a few remarks.
\br
The expression \eqref{ldf} coincides with the prediction of Macroscopic
Fluctuation Theory with transport coefficients
$$
D(\rho) = \frac{1}{2s}, \qquad  \qquad \sigma(\rho)= \frac{\rho}{2s} \left(1+\frac{\rho}{2s}\right)
$$
which indeed are the transport coefficient of the harmonic model, as proved in
\cite{capanna}. In particular, for $s=1/2$, we recover the transport coefficient of
the discrete KMP model and the large deviation function  \eqref{ldf}
coincides with the one computed in \cite{bgl}.
There it was already remarked that the infimum over $\theta$ can be viewed as
a contraction principle over a random local temperature profile given by uniform
order statistics. The macroscopic fluctuation theory can strictly speaking not be
applied to the KMP model, or to any of the models studied in this paper, because
the proof requires superexponential tails of the marginals of the equilibrium product
measures, which does not hold for any of the models in the KMP class.
Therefore, even if Theorem \ref{theo-ld} gives the large deviation principle for the
whole class of harmonic models with parameter $2s$ integer, it does not prove yet the same for the KMP
model and its generalizations. Nevertheless  the macroscopic fluctuation theory predicts
that these models sharing the same macroscopic transport coefficients have the same
rate function.
\er
\br
As already remarked in \cite{bgl} for the case $s=1/2$, the rate function  \eqref{ldf} is
non-convex.  This is at the root of the fact that the large deviation function
can not be represented as the Legendre transform of a convex function.
Indeed if one takes the Legendre transform of the pressure
one rather obtains the convex hull of the rate function.
\er
\br
For the models with compact state space, such as the exclusion process,
the expression for the large deviation function contains a supremum,
rather then an infimum \cite{DLS, mft}. For the weakly asymmetric exclusion
process the density large deviation has been written as a minimization problem
(see formula (2.3) of \cite{ed}) and for the asymmetric exclusion process a
contraction involving Brownian excursions has been considered \cite{del}.
\er
\bpr[Proof of Theorem \ref{theo-ld}]
Preliminarily, consider an inhomogenous product measure with marginal Negative Binomials with a smooth slowly varying parameter.
Thus, assume we have a measure $\mu_N$ of the form
\be\label{wiwop}
\mu_N = \otimes_{i=1}^N \nu_{\theta(\frac{i}{N})}
\ee
where $\nu_{\theta(\frac{i}{N})}$ is the Negative Binomial measure introduced in \eqref{negbin-ew} with mean $\theta(\frac{i}{N})$
and where $\theta:[0,1]\to[0,\infty)$ is a smooth increasing function.
We call
\be\label{empprof}
\ell_N= \frac1N\sum_{i=1}^N\eta_i \delta_{i/N}
\ee
the empirical density profile when $\eta$ has distribution $\mu_N$.
Then, G\"artner-Ellis theorem tells us that the sequence of measures $(\ell_N)_{N\ge 1}$
satisfies a large deviation principle with a good rate function
${I}(\rho,\theta)$. The LDP of $(\ell_N)_{N\ge 1}$ has to be interpreted
 in the set of positive finite measures on $[0,1]$ equipped
with the weak topology. We have ${I}(\rho,\theta)=\infty$ for a measure $\rho$
which is not absolutely continuous w.r.t. Lebesgue measure on $[0,1]$;
otherwise the rate function $I(\rho,\theta)$ is given
and is obtained as the Legendre transform of the pressure
\beq\label{progom}
{I}(\rho,\theta)= \sup_{h}\left( \int \rho(x) h(x) dx- P(h,\theta)\right)
\eeq
where
\beq
P(h,\theta)&=& \lim_{N\to\infty} \frac1N\log\E_{\mu_N}\left(e^{N\langle \ell_N, h\rangle}\right)
\nn\\
&=& \lim_{N\to\infty} \frac1N\log\E_{\mu_N} e^{\sum_{i=1}^N \eta_x h(i/N)}\nonumber\\
\eeq
has been computed in \eqref{pressure-h-theta}. Evaluating the Legendre transform
one obtains for ${I}(\rho,\theta)$ the expression that is given in \eqref{i-rho-theta}.

The type of measures which are of interest to us, are not product measures of the form \eqref{wiwop},
but product measures with  parameters  that are themselves random variables.
More precisely we have a measure of the form
\be\label{wiwi}
\mu_{N,\rho_l,\rho_r}= \mathbb{E}\left(\otimes_{i=1}^N \nu_{\Theta_{2si,n}}\right)
\ee
where  $n=2s(N+1)-1$ and the additional expectation refers to the random variables $\Theta_{1,n} \le \Theta_{2,n}  \le \ldots \le \Theta_{n,n}$
which are the ascending order statistics of a sequence $\Theta_1\ldots,\Theta_n$
of  i.i.d. random variables with common uniform  distribution on the interval $[\rho_l,\rho_r]$.
Recalling  the definition of the sample path of the order statistics
$$
\Theta_{n}(x) = \Theta_{\lfloor (n+1)x\rfloor +1,n} \qquad x\in[0,1], \qquad \qquad\text{with} \quad \Theta_{n+1,n}:=\rho_r$$
the stationary measure is rewritten as
\be\label{wiwi2}
\mu_{N,\rho_l,\rho_r}= \mathbb{E}\left(\otimes_{i=1}^N \nu_{\Theta_{n} \left(  \frac{i}{N+1} \right) } \right) \;.
\ee
As we know from Lemma \ref{ld-order},  the sample path of the order statistics of uniform random variables satisfies a large deviation principle
with rate function $J(\theta)$ given in \eqref{Jt}.  As a consequence, the contraction principle gives that,
under $\mu_{N,\rho_l,\rho_r}$, the sequence $(L_N)_{N\ge 1}$
satisfies the large deviation principle with rate function $I$ which is
only finite on positive measures $\rho$ of the form $\rho(x) dx$,
where it is equal to
$$
 {I}(\rho) = \inf_{\underset{\theta(1)=\rho_r}{\underset{\theta(0) = \rho_l}{\underset{\text{strictly increasing }}{\theta:[0,1]\to \mathbb{R}_+}}}}
 \Big[I(\rho,\theta) + J(\theta)\Big]\, .
$$
\epr

\section{Additivity Principle}

In this section we compare the moment generating function of system of size $N$
to the moment generating function of two subsystems of sizes $N_1, N_2$ with
$N_1+N_2 = N$.  In the macroscopic limit (i.e.~when the two subsystems are of macroscopic sizes
$N_1 = N x$ and $N_2 = N(1-x)$ with $x\in(0,1)$) we get a rigorous proof of an additivity
principle for the pressure (and similarly for the density large deviations).
In the non-equilibrium set-up, an additivity principle was
first established in \cite{DLS} for the density profile large deviations of the non-equilibrium
steady state of the symmetric exclusion process. Surprisingly, the corresponding
additivity principle for the pressure of the symmetric exclusion process contained an infimum,
whose physical basis remain not understood. The pressure additivity principle proved
here for the harmonic model contains instead a supremum and generalizes the one
conjectured in \cite{bgl} for the discrete-KMP model. The proof relies on
an integral equation  (see \eqref{integral-eq2} below)
relating the partition functions of the systems of
sizes $N_1,N_2$ and $N$ and an application of Varadhan's lemma.
The integral equation is in turn a consequence of the properties of order statistics,
in particular the Markovian structure of Lemma \ref{lemma2-order} and the properties
of conditioning of Lemma \ref{lemma3-order} and Lemma \ref{lemma4-order}.

As it will be discussed in Section \ref{section-constant} , the additivity principle for the pressure implies that the pressure for constant field, corresponding to the large deviations of the total density, 
determines completely the pressure of any other field, by approximation by piece-wise constant functions. 
This implies in particular that Theorem \ref{pres1thm} completely determines the pressure.

\subsection{The additivity principle for the pressure}

In order to formulate the additivity principle, we need to generalise the definition of pressure given in \eqref{PP} to the case of a system whose macroscopic volume is the interval $[a,b]$ and the
boundary densities are $0<\rho_a\le \rho_b$. This is obtained by starting from a microscopic system with $\lceil (b-a)N\rceil$ sites and taking the limit as $N\to \infty$
\be\label{PPP}
{P}^{[a,b]}_{\rho_a,\rho_b}(h) :=  \lim_{N\to\infty} \frac{1}{ N} \log \mathbb{E} \Big[e^{\sum_{i=1}^{N_{a,b}} \eta_i h\left(a+\tfrac{i}{N}\right)} \Big] \, \qquad \text{with} \qquad N_{a,b}= \lceil (b-a)N\rceil\,.
\ee
Here $h:[a,b]\to \mathbb{R}$ and $\mathbb{E}$ denotes  expectation with respect to the stationary measure $\mu_{N_{a,b}, \rho_a, \rho_b}$.
As we did in  \eqref{EEE} for the system with macroscopic unit volume, the expectation in \eqref{PPP} can be written in terms of the moment generating function:
\be
\label{mix-2s-intab}
\Psi^{[a,b]}_{N_{a,b},\rho_a,\rho_b}(\mathsf{h})
 = \mathbb{E}\bigg[\prod_{i=1}^{N_{a,b}} M_{\Theta_{2si,n_{a,b}}}(\mathsf{h}_i)\bigg], \qquad \text{with} \qquad n_{a,b}= 2s(N_{a,b}+1)-1
\ee
defined on vectors $\mathsf{h} \in {\cal A}_{N_{a,b},\rho_a,\rho_b}$. Here $\Theta_{1,n_{a,b}}\le \Theta_{2,n_{a,b}}\le  \ldots \le \Theta_{n_{a,b},n_{a,b}}$ is the ascending order statistics of $n_{a,b}$ independent uniform random variables on $[\rho_a, \rho_b]$ and and $M_{\theta}(\cdot)$ is the moment generating function of a Negative
Binomial distribution with parameters $(2s,\theta)$, as defined in \eqref{pgf-negbin}.
It then follows that
\be\label{ppp}
{P}^{[a,b]}_{\rho_a,\rho_b}(h) :=   \lim_{N\to\infty} \frac{1}{ N} \log \Psi^{[a,b]}_{N_{a,b},\rho_a,\rho_b}({\mathsf{h}}^{(N)})
\ee
where $\mathsf{h}^{(N)}$ is the $N_{a,b}$-dimensional  vector of components:
\be
\mathsf{h}^{(N)}_i:= h\left(a+\tfrac iN \right), \qquad \text{for} \qquad i=1, \ldots N_{a,b} \;.
\ee
Furthermore, to formulate the additivity principle, we define the modified pressure
\be
\widetilde P^{[a,b]}_{\rho_a,\rho_b}(h) := P^{[a,b]}_{\rho_a,\rho_b}(h) + 2s (b-a) \log\Big(\frac{\rho_b-\rho_a}{b-a}\Big)\,.
\ee
In the next theorem we prove that the modified pressure satisfies an the additivity principle.
\bt[Pressure additivity principle]
\label{adprop}
Let $0<\rho_l < \rho_r $, $0<x<1$ and $h: [0,1]\to \mathbb R$, then we have
\be\label{adprn2}
\widetilde P^{[0,1]}_{\rho_l,\rho_r}(h)
=\sup_{\rho_l\leq \theta\leq \rho_r} \Big[\widetilde P^{[0,x]}_{\rho_l,\theta}(h_1)  + \widetilde P^{[x,1]}_{\theta,\rho_r}(h_2) \Big]
\ee
where $h_1: [0,x]\to \mathbb R$ and $h_2: [x,1]\to \mathbb R$ are the restrictions of $h$ respectively, to $[0,x]$ and to $[x,1]$.
More generally,
for  $\kappa \ge 2$ and $0=x_0 \le x_1 \le \ldots \le x_{\kappa} = 1$,  calling $h_i: [x_{i-1},x_i]\to \mathbb R$ the restriction of $h$ to $[x_{i-1},x_i]$,  for $i=1, \ldots, \kappa$, we have
\be\label{adgeneral}
\widetilde P^{[0,1]}_{\rho_l,\rho_r}(h)
=\sup_{\rho_0\leq \rho_1\leq\ldots\leq\rho_{\kappa-1} \le \rho_{\kappa}} \sum_{i=1}^{\kappa} \widetilde P^{[x_{i-1},x_{i}]}_{\rho_{i-1}, \rho_i}(h_i)
\ee
with the convention $\rho_0=\rho_l$, $\rho_{\kappa}=\rho_r$.

\et

\bpr
We prove \eqref{adprn2}, i.e. the case $\kappa=2$, the  case  of a generic $\kappa$ can be then deduced by induction.
\noindent
As a first step we fix two integers $N_1,N_2 \in \mathbb{N}$
 such that $N_1+N_2=N$ and prove the following identity for the moment generating function
\be\label{integral-eq2}
\Psi_{N,\rho_l,\rho_r}(\mathsf{h}_1,\ldots,\mathsf{h}_N) =
\mathbb{E}\Big(M_{\Theta_{2sN_1,n_1}}({\mathsf{h}_{N_1}})
\Psi_{N_1-1,\rho_l,\Theta_{2sN_1,n_1}}(\mathsf{h}_1,\ldots, \mathsf{h}_{N_1-1})
\Psi_{N_2,\Theta_{2sN_1,n_2},\rho_r}(\mathsf{h}_{N_1+1},\ldots,\mathsf{h}_{N})
\Big)
\ee
where $n_1= 2sN_1 -1$, $n_2= 2s(N_2+1) -1$. Here  $\Theta_{2sN,n}$ is the $2s{N}^{\text{th}}$-th ascending order statistics of $n$ independent uniforms
on the interval $(\rho_l,\rho_r)$.
\vskip.1cm
\noindent
In order to prove \eqref{integral-eq2} we start from Proposition \ref{prop-conclude}  which says that, for $\mathsf{h} \in {\cal A}_{N,\rho_l,\rho_r}$,
$$
\Psi_{N,\rho_l,\rho_r}(\mathsf{h}) = \mathbb{E}\left[ \prod_{i=1}^N M_{\Theta_{2si,n}}(\mathsf{h}_i)\right]
$$
with
$$
\Theta_{2si,n} = \rho_l + (\rho_r-\rho_l) U_{2si,n}, \qquad \quad i=1, \ldots, n
$$
where $U_{2si,n}$ is the $2si^{\text{th}}$ order statistics
of $n=2s(N+1) -1$ i.i.d. random variables that are uniformly distributed on the interval $(0,1)$.
The tower property of conditional expectation implies
\begin{eqnarray}
\Psi_{N,\rho_l,\rho_r}(\mathsf{h}_1,\ldots, \mathsf{h}_N)
& = &
 \mathbb{E}\left(\mathbb{E}\left(  \prod_{i=1}^N  M_{\Theta_{2si,n} }(\mathsf{h}_i) \,\Big|\, \Theta_{2sN_1,n} \right) \right) \nn \\
& = &
\label{primo}
\mathbb{E}\left(M_{\Theta_{2sN_1,n}}(\mathsf{h}_{N_1}) \, \mathbb{E}\left(  \prod_{\substack{i =1 \\ i \neq N_1}}^N  M_{\Theta_{2si,n} }(\mathsf{h}_i)\,\Big|\, \Theta_{2sN_1,n} \right)\right) \;.
\end{eqnarray}
Now, given $\theta\in[\rho_l,\rho_r]$, the event $\{ \Theta_{2sN_1,n} = \theta\}$ is equivalent
to the event $\{ U_{2s N_1,n} = u\}$, with $u= \frac{\theta-\rho_l}{\rho_r - \rho_l}$.
Therefore, using  the property  of the conditional expectations of
order statistics (specifically equation \eqref{main-order} of Lemma \ref{lemma4-order}
 with $n=2s(N+1) -1$ and $m=2sN_1$) we obtain
\begin{eqnarray}
\label{secondo}
\mathbb{E}\left(  \prod_{\substack{i =1 \\ i \neq N_1}}^N  M_{\Theta_{2si,n}}(\mathsf{h}_i)\,\Big|\, \Theta_{2sN_1,n} =\theta\right)
& = &
\mathbb{E}\left(  \prod_{\substack{i =1 \\ i \neq N_1}}^N  M_{\Theta_{2si,n}}(\mathsf{h}_i)\ \,\Big|\, U_{2sN_1,n} = u \right)\nn \\
& = &
\mathbb{E}\left(  \prod_{i=1}^{N_1-1}  M_{\Theta^{\star}_{2si,n_1}}(\mathsf{h}_i)  \right) \cdot
\mathbb{E}\left(  \prod_{i=1}^{N_2}  M_{\tilde{ \Theta}_{2si,n_2}}(\mathsf{h}_{N_1+i}) \right)
\end{eqnarray}
where
$$
\Theta^{\star}_{2si,n_1} = \rho_l + (\rho_r-\rho_l) U^{\star}_{2si,n_1} \qquad\qquad i=1,\ldots,N_1-1
$$
with $U^{\star}_{2si,n_1}$ the $2si^{\text{th}}$ order statistics
of $n_1= 2sN_1 -1$ i.i.d. random variables uniformly distributed on the interval $(0,u)$
and similarly
$$
\tilde{ \Theta}_{2si,n_2}= \rho_l + (\rho_r-\rho_l) \tilde{ U}_{2si,n_2} \qquad\qquad i=1,\ldots,N_2
$$
with $\tilde{U}_{2si,n_2}$ the $2si^{\text{th}}$ order statistics
of  $n_2 = 2s(N_2+1)-1$ i.i.d. random variables that are uniformly distributed on the interval $(u,1)$.
In other words, defining
$$
\theta(u) = \rho_l + u(\rho_r-\rho_l),\qquad u\in [0,1]
$$
 the $\{\Theta^{\star}_{2si,n_1}\}_{i=1,\ldots,N_1-1}$ are the order statistics
 (sampled every $2s$ steps) of $n_1= 2sN_1 -1$  i.i.d. uniforms on $(\rho_l,\theta(u))$
and  the $\{\tilde{ \Theta}_{2si,n_2}\}_{i=1,\ldots,N_2}$ are the order statistics
(sampled every $2s$ steps) of $n_2= 2s(N_2+1) -1$  i.i.d. uniforms on $(\theta(u),\rho_r)$.
As a consequence, combining \eqref{primo} and \eqref{secondo}, we obtain \eqref{integral-eq2}.
\vskip.1cm
\noindent
We further proceed by observing that, recalling    \eqref{pgf-negbin}, the identity \eqref{integral-eq2} can be explicitly written as
\begin{eqnarray}\label{integral-eq333}
\Psi_{N,\rho_l,\rho_r}(\mathsf{h}_1,\ldots, \mathsf{h}_N)
&= &
\int_{\rho_l}^{\rho_r} d\theta \:
\Psi_{N_1-1,\rho_l,\theta}(\mathsf{h}_1,\ldots, \mathsf{h}_{N_1-1})\cdot
\Psi_{N_2,\theta,\rho_r}(\mathsf{h}_{N_1+1},\ldots, \mathsf{h}_{N})\nn\\
& & \qquad \cdot \bigg( \tfrac{1}{1+(1-e^{\mathsf{h}_{N_1}})\theta}\bigg)^{2s} \cdot \tfrac{1}{\rho_r-\rho_l}\cdot f_{U_{2sN_1,n}}\left(\tfrac{\rho_r-\theta}{\rho_r-\rho_l}\right)
\end{eqnarray}
where $f_{U_{2sN_1,n}}$ is the probability density of the random variable $U_{2sN_1,n}$ which, from Lemma  \ref{lemma1-order}, is equal to
\be
f_{U_{2sN_1,n}}(u)= \frac{(2s(N+1)-1)!}{(2s N_1-1)!(2s(N_2+1) -1)!} \cdot u^{2sN_1-1} (1-u)^{2s(N_2+1) -1} \,.
\ee
In order to take the macroscopic limit we  consider blocks of macroscopic sizes i.e. $N_1 = \lfloor N x \rfloor$ and $N_2= \lfloor N(1-x)  \rfloor$, with $x\in (0,1)$. Now  let $h:[0,1]\to \mathbb R$ and let $h_1: [0,x]\to \mathbb R$ and $h_2: [x,1]\to \mathbb R$ be the restrictions of $h$ to $[0,x]$ and to $[x,1]$.
Then by definition we have that
\be
 \lim_{N\to\infty} \frac{1}{ N} \log \Psi^{[0,x]}_{ \lfloor N x \rfloor-1,\rho_l,\theta}\left(h\left(\tfrac 1 N\right), \ldots, h\left(\tfrac { \lfloor N x\rfloor -1} {N}\right) \right)={P}^{[0,x]}_{\rho_l,\theta}(h_1) \nonumber
\ee
and
\be
 \lim_{N\to\infty} \frac{1}{N} \log \Psi^{[x,1]}_{ \lfloor N (1-x) \rfloor,\theta,\rho_r}\left(h\left(\tfrac { \lfloor N (1-x)\rfloor } {N}\right), \ldots, h\left(\tfrac{N}{N})\right) \right)=
 {P}^{[x,1]}_{\theta,\rho_r}(h_2) \, .\nonumber
\ee
Moreover, using that
\beq
f_{U_{2s \lfloor N x \rfloor,n}}(u)= e^{2sN[x\log \tfrac ux +(1-x)\log \tfrac{1-u}{1-x} +o(1)]}\nonumber
\eeq
and considering \eqref{integral-eq333} for a vector with components $\mathsf{h}^{(N)}_i:= h\left(a+\tfrac iN \right)$ with  $i=1, \ldots N_{a,b}$
we obtain
\beq
\label{integral-eq3}
\Psi^{[0,1]}_{N,\rho_l,\rho_r}\left(h\left(\tfrac 1 N\right), \ldots, h\left(\tfrac {N} {N}\right) \right)
 =
\int_{\rho_l}^{\rho_r}
e^{N \left[ P^{[0,x]}_{\rho_l,\theta}(h_1)  +{P}^{[x,1]}_{\theta,\rho_r}(h_2) + 2s x\log \frac{\theta-\rho_l}{x(\rho_r - \rho_l)} +2s (1-x)\log \frac{\rho_r-\theta}{(1-x)(\rho_r - \rho_l)}+o(1) \right]} \; d\theta\nonumber\,.
\eeq
Then, taking the limit as $N\to \infty$ and recalling the definition of the modified pressure, the claim \eqref{adprn2} follows from the Laplace principle.

\epr

\subsection{The additivity principle for the density large deviation function}
\label{section-large}

For a macroscopic system  on the interval  $[a,b]$ we define the modified density large deviation function
with boundary parameters $0<\rho_a < \rho_b$ as
\be
\widetilde I^{[a,b]}_{\rho_a,\rho_b}(\rho) := I^{[a,b]}_{\rho_a,\rho_b}(\rho) - 2s (b-a) \log\Big(\frac{\rho_b-\rho_a}{b-a}\Big)
\ee
where $I^{[a,b]}_{\rho_a,\rho_b}(\cdot)$ is the large deviation function of the empirical profile
$$
L_N^{[a,b]} = \frac{1}{N_{a,b}} \sum_{i=1}^{N_{a,b}} \eta_i \delta_{a +\frac{i}{N}}.
$$
\bt[Large deviation additivity principle]
\label{adprop22}
For  $0<x<1$ and $\rho: [0,1]\to \mathbb R$,
we have
\be\label{adprld}
\widetilde I^{[0,1]}_{\rho_l,\rho_r}(\rho)
=\inf_{\rho_l\leq \theta\leq \rho_r} \Big[\widetilde I^{[0,x]}_{\rho_l,\theta}(\rho_1)  + \widetilde I^{[x,1]}_{\theta,\rho_r}(\rho_2) \Big]
\ee
where $\rho_1: [0,x]\to \mathbb R$ and $\rho_2: [x,1]\to \mathbb R$ are the restrictions of $\rho$ respectively, to $[0,x]$ and to $[x,1]$.
More generally, for $\kappa \ge 2$ and
for  $0=x_0 \le x_1 \le \ldots \le x_\kappa = 1$, calling $\rho_i: [x_{i-1},x_i]\to \mathbb R$ the restriction of $\rho$ to $[x_{i-1},x_i]$,  for $i=1, \ldots, \kappa$, we have
\be\label{adgeneral22}
\widetilde I^{[0,1]}_{\rho_l,\rho_r}(\rho)
=\inf_{\theta_0\leq \theta_1\leq\ldots\leq\theta_{\kappa-1} \le \rho_{\kappa}} \sum_{i=1}^{\kappa} \widetilde I^{[x_{i-1},x_{i}]}_{\theta_{i-1}, \theta_i}(\rho_i)
\ee
with the convention $\theta_0=\rho_l$, $\theta_{\kappa}=\rho_r$.
\bpr
The proof is analogous to the  one of Theorem \ref{adprop}.
\epr
\et

\section{Explicit formulas for the pressure and further results on the additivity principle}
\label{section-constant}

In this final section, we give explicit formulas for the pressure and prove equivalence between the additivity principle and
the MFT variational expression.
Firstly, in subsection \ref{section-constant1} we find an explicit formula for the pressure in a constant field. In the spirit of this paper, we show
how this can be achieved in two ways: either macroscopically, solving the MFT variational principle, or
microscopically, using the explicit characterization of the stationary measure to produce upper and lower bounds matching in the limit
$N\to\infty$. Secondly, in subsection \ref{section-constant2}, using the knowledge of the pressure in a constant field,
we prove the equivalence between Theorem 4.1 (pressure MFT variational problem) and Theorem 6.1 (pressure additivity principle).
Thirdly, in subsection \ref{section-constant3},
we consider the finite-volume pressure $P_N$ for a constant field. We prove that it
satisfies a recursion relation in $N$, which in fact can be solved for the Laplace transform.
In particular, we prove that the finite-volume pressure of the model with $s=1/2$ is
size-independent, i.e. it takes the same value for all system sizes $N$.

\subsection{The pressure for a constant field}
\label{section-constant1}

We analyse in detail the case of constant  field, i.e. $h(x)= \mathsf{h}\in \mathbb R$ for all $x\in [0,1]$.

\subsubsection{Solution of MFT variational problem}

When the field $h(\cdot)$ is constantly equal to $\mathsf{h}$, the variational problem for the pressure reads
\be \label{pressuresup}
P^{[0,1]}_{\rho_l,\rho_r} (\mathsf{h}) = \sup_{\theta} {\cal P}  (\mathsf{h}, \theta)
\ee
with
\be
\label{pressureMFT-fct0}
 {\cal P} (\mathsf{h}, \theta) = 2s \int_{0}^1 dx \left[
 \log \Big(
 \frac{1}{1+(1-e^{\mathsf{h}})\theta(x)} \Big)
+  \log\Big(
 \frac{\theta'(x)}{\rho_r -\rho_l}\Big) \right]
\ee
and the supremum is over all functions $\theta:[0,1]\to \mathbb{R}$ monotone such that $\theta(0)= \rho_l$ and $\theta(1)=  \rho_r$.
In other words
$$
\label{press-optimal}
P^{[0,1]}_{\rho_l,\rho_r} (\mathsf{h}) =  {\cal P} (\mathsf{h}, \theta_*)
$$
where $\theta_*$ is defined implicitly by $\frac{\delta {\cal P}}{\delta \theta}\Big|_{\theta=\theta_*} = 0$.
Computing the functional derivatives one gets the boundary value problem
\be
\label{eqdiff2}
\frac{1 - e^\mathsf{h}}{1+(1-e^{\mathsf{h}})\theta_*}  - \frac{ \theta''_*}{(\theta'_*)^2} =0, \qquad \theta_*(0)= \rho_l, \: \theta_*(1)=   \rho_r
\ee
whose solution is given by
\be\label{thetastar}
\theta_*(x)=  \frac{1}{1-e^\mathsf{h}} \Big[ (\rho_l (1-e^\mathsf{h}) + 1) \left(\frac{\rho_r (1-e^\mathsf{h}) + 1}{\rho_l (1-e^\mathsf{h}) + 1} \right)^x -1\Big].
\ee
Plugging \eqref{thetastar} in \eqref{pressureMFT-fct0} one obtains
\be
\label{press}
P^{[0,1]}_{\rho_l,\rho_r} (\mathsf{h}) = {\cal P} (\mathsf{h}, \theta_*) = 2s \log\left(\frac{1}{(\rho_r-\rho_l)(1-e^\mathsf{h})} \log \frac{1+(1-e^\mathsf{h})\rho_r}{1+(1-e^\mathsf{h})\rho_l}\right).
\ee
In a similar manner, it can be proved that
\be
\label{pressab}
P^{[a,b]}_{\rho_a,\rho_b} (\mathsf{h}) =2s(b-a)\log \left( \frac 1{(\rho_b-\rho_a)\left(1-e^h\right)} \cdot  \log\frac{1+\rho_b(1-e^\mathsf{h})}{1+\rho_a(1-e^\mathsf{h})}\right)\nn
\ee
where $P^{[a,b]}_{\rho_a,\rho_b}(\cdot)$ is the pressure for a  system in the macroscopic interval $[a,b]$.

\subsubsection{Matching upper and lower bound}

In this section we consider the
 moment generating function evaluated in a point with components all equal to each others, i.e. $(\mathsf{h},\ldots, \mathsf{h})$, with $\mathsf{h} \in \mathbb{R}$ . For this observable  we introduce the  notation
 $\Psi_{N,\rho_l,\rho_r}^{(1)}: \mathbb{R} \to \mathbb{R}$ for the one-variable function
 \be
 \Psi_{N,\rho_l,\rho_r}^{(1)}(\mathsf{h}):= \Psi_{N,\rho_l,\rho_r}(\mathsf{h}, \ldots, \mathsf{h}).
 \ee
From Proposition \ref{prop-conclude} we know that, thanks to the mixture structure of the non-equilibrium steady state, this can be written as
\begin{equation*}
\Psi^{(1)}_{N,\rho_l,\rho_r}(\mathsf{h}) = \mathbb{E} \left[ \prod_{i=1}^{N} M^{2s}_{\Theta_{2si,n}}(\mathsf{h}) \right]
\end{equation*}
where we recall that, for $2s \in \mathbb N$, $M^{2s}_{\theta}(\cdot)$ is the generating function of a Negative Binomial of parameters $2s$ and $\theta$, i.e.
\be
M^{1}_{\theta}(\mathsf{h})= \frac 1 {1+\theta (1-e^{\mathsf{h}})} \qquad \text{and} \qquad M^{2s}_{\theta}(\mathsf{h})= (M^{1}_{\theta}(\mathsf{h}))^{2s}.
\ee
Notice that we added the superscript $2s$ in the notation for this generating function because in what follows it will be crucial to distinguish the case of general $2s \neq 1$ and $2s = 1$.
In the following theorem we will prove that the logarithm of $\Psi_{N,\rho_l,\rho_r}^{(1)}(\mathsf{h})$ divided by $N$ converges, in the limit as $N\to \infty$, to the solution of the variational problem for the pressure given in \eqref{press}. We will restrict to the case $2s \in \mathbb N$.

\bt[Pressure, constant field]
\label{pres1thm}
For all  $s>0$ with $2s \in \mathbb{N}$, $\mathsf{h} \in \mathbb R$ we have that
\be\label{pres1}
\lim_{N\to \infty} \frac{1}{N}\log \Psi^{(1)}_{N,\rho_l,\rho_r}(\mathsf{h}) = 2s\log\left(\frac{1}{(\rho_r - \rho_l)(1-e^{\mathsf{h}})}\log\frac{1+(1-e^{\mathsf{h}})\rho_r}{1+(1-e^{\mathsf{h}})\rho_l}\right)=P^{[0,1]}_{\rho_l,\rho_r} (\mathsf{h})  \;.
\ee
\et
\bpr
Consider first $2s=1$. In this case, because $n := 2s(N+1)-1=N$, the joint distribution
of $(U_{2s,n},\ldots, U_{2sN,n})$ is simply the joint distribution of the order statistics
$(U_{1,N}, \ldots, U_{N,N})$. As a consequence, the corresponding variables $\Theta_{1,N} \ldots, \Theta_{N,N}$ defined in \eqref{rho-i} are
the order statistics of $N$ uniforms on the interval $[\rho_l,\rho_r]$.
Let us consider  $N$ independent uniform random variables on the interval $[\rho_l,\rho_r]$, denoted $\Theta_1, \ldots, \Theta_N$ as in equation \eqref{rho-i}. Then for every smooth function
$g$ we have that in distribution,
\[
\prod_{i=1}^N g(\Theta_{i,N})= \prod_{i=1}^N g(\Theta_i)
\]
because in the product of all the $N$ terms the ordering does not matter.
As a consequence,
\[
 \Psi^{(1)}_{N,\rho_l,\rho_r}(\mathsf{h}) =  \mathbb{E} \left(\prod_{i=1}^N M^{1}_{\Theta_{i,N}}(\mathsf{h})\right)=\mathbb{E} \left(\prod_{i=1}^N M^{1}_{\Theta_i}(\mathsf{h}) \right)= \Big[ \mathbb{E}\left( M^{1}_{\Theta_1}\left( \mathsf{h} \right) \right) \Big] ^{N}
\]
where in the last step we used independence of the $\Theta_i$.
Since
\begin{equation} \label{t1}
\E (M^{1}_{\Theta_1}(\mathsf{h}))= \frac{1}{\rho_r-\rho_l}\int_{\rho_l}^{\rho_r} \frac{d\rho}{1+\rho(1-e^{\mathsf{h}})}= \frac{1}{(\rho_r - \rho_l) (1-e^{\mathsf{h}})}\log\frac{1+(1-e^{\mathsf{h}})\rho_r}{1+ (1-e^{\mathsf{h}}) \rho_l} \;,
\end{equation}
we immediately get the result for the infinite pressure
\begin{equation*}
\lim_{N\to \infty} \frac{1}{N}\log \Psi^{(1)}_{N,\rho_l,\rho_r}(\mathsf{h}) = \log\left(\frac{1}{(\rho_r - \rho_l)(1-e^{\mathsf{h}})}\log\frac{1+(1-e^{\mathsf{h}})\rho_r}{1+(1-e^{\mathsf{h}})\rho_l} \right)  \;.
\end{equation*}
To deal with the general case, first notice that
the joint distribution of $(\Theta_{2s,n}, \ldots \Theta_{2sN,n})$ can be
obtained as follows. We consider $n:=2s(N+1)-1$ independent uniforms $(\Theta_1, \ldots, \Theta_{n})$ on the interval $[\rho_l, \rho_r]$ and denote by $(\Theta_{1,n}, \ldots, \Theta_{n,n})$ the ordered vector.
By sampling every $2s$ steps implies that $(\Theta_{2s,n}, \ldots, \Theta_{2sN,n})$ is equal in distribution to $(\Theta_{1,N}, \ldots \Theta_{N,N})$. Moreover, $M^{2s}_{\theta}({\mathsf{h}})= \left( M^{1}_{\theta}(\mathsf{h}) \right) ^{2s}$
and therefore
\[
\prod_{i=1}^N M^{2s}_{\Theta_{i,n}}(\mathsf{h})= \prod_{i=1}^N \left(  M^{1}_{\Theta_{2si,n}}(\mathsf{h}) \right) ^{2s} \;.
\]
We notice that for $\mathsf{h}$ fixed, the function $\theta \to M^{1}_{\theta}(\mathsf{h})$ is non-decreasing and bounded
from above and below by positive constants, i.e.,
$$
0<c_1\leq M^{1}_{\theta}(\mathsf{h}) \leq c_2<\infty \;.
$$
As a consequence,
\begin{equation*}
 \Psi^{(1)}_{N,\rho_l,\rho_r}(\mathsf{h}) = \mathbb{E}\left( \prod_{i=1}^{N} M^{2s}_{\Theta_{i,N}} (\mathsf{h}) \right)  = \mathbb{E}\left( \prod_{i=1}^{N} \left(  M^{1}_{\Theta_{2si,n}} (\mathsf{h})\right)^{2s}  \right) \geq \mathbb{E}\left( \prod_{i=1}^{2sN}  M^{1}_{\Theta_{i,n}} (\mathsf{h})  \right)
\end{equation*}
where the last inequality follows from the fact that  $M^{1}_{\cdot} (\mathsf{h}) $ is non-decreasing and  for $i=1, \ldots, N$
$$
\Theta_{2si,n} \ge \Theta_{j,n} \qquad \text{when} \qquad 2s(i-1) < j \leq 2si \;.
$$
Considering the $\log$, dividing by $N$ and taking the $N \to \infty$ limit on both sides, we have
\begin{equation*}
\lim_{N\rightarrow \infty} \dfrac{1}{N} \log  \left(  \Psi^{(1)}_{N,\rho_l,\rho_r}(\mathsf{h})\right)  \geq
\lim_{N\rightarrow \infty} \dfrac{1}{N} \log  \mathbb{E}\left( \prod_{i=1}^{2sN} \left(  M^{1}_{\Theta_{i,n}} (\mathsf{h})\right)^{2s}  \right) = \lim_{N\rightarrow \infty} \dfrac{1}{N} \log  \mathbb{E}\left( \prod_{i=1}^{2s(N+1) -1} \left(  M^{1}_{\Theta_{i,n}} (\mathsf{h})\right)^{2s}  \right)
\end{equation*}
where the last identity follows from the boundedness of $M^{1}_.$, which is used to add $2s-1$ terms in the product.
As for $s=1/2$ we can now remove the order and use the independence of the $\Theta_i$, $i=1, \ldots, n$
\begin{align*}
\lim_{N\rightarrow \infty} \dfrac{1}{N} \log  \left(  \Psi^{(1)}_{N,\rho_l,\rho_r}(\mathsf{h})\right) & \geq  \lim_{N\rightarrow \infty} \dfrac{1}{N} \log  \mathbb{E}\left( \prod_{i=1}^{2s(N+1) -1} \left(  M^{1}_{\Theta_{i,n}} (\mathsf{h})\right)^{2s}  \right)  \\ & =
\lim_{N\rightarrow \infty} \dfrac{1}{N} \log  \mathbb{E}\left( \prod_{i=1}^{2s(N+1) -1} \left(  M^{1}_{\Theta_{i}} (\mathsf{h})\right)^{2s}  \right) \\ & = \lim_{N\rightarrow \infty} \dfrac{1}{N} \log \left(  \mathbb{E} \left( M^{1}_{\Theta_1} (\mathsf{h})  \right) \right)^{2s(N+1)-1}
 \\ & = 2s \log \left( \frac{1}{(\rho_r - \rho_l)(1-e^{\mathsf{h}})}\log\frac{1+(1-e^{\mathsf{h}})\rho_r}{1+(1-e^{\mathsf{h}})\rho_l} \right)
\end{align*}
where the last identity follows from \eqref{t1}.

The idea to obtain a matching upper bound is similar.
Now, for $i = 1, \ldots, N$, we consider $2si \leq j \le 2s(i+1)$  so that $ \Theta_{2si,n} \leq  \Theta_{j,n}$ implies
\begin{equation*}
\mathbb{E} \left( \prod_{i=1}^{N} \left( M^{1}_{ \Theta_{2si,n}} (\mathsf{h}) \right) ^{2s} \right) \leq
\mathbb{E} \left( \prod_{i=2s}^{2s(N+1) -1} M^{1}_{\Theta_{i,n}} (\mathsf{h})  \right) \;.
\end{equation*}
since $M^{1}_{\cdot} (\mathsf{h}) $ is non-decreasing.
As before, in the limit we can consider the full product from $i=1, \ldots, 2s(N+1) -1$ by adding the first $2s-1$ terms so that we can replace the ordered variables $\Theta_{i,n}$ with the corresponding non ordered ones $\Theta_{i}$ and use their independence to conclude the proof.
\epr

\br[Case $s=1/2$]
In the course of the previous proof, we have proven, in particular that, for the case $s=1/2$ the  constant field generating function $\Psi^{(1)}_{N,\rho_l,\rho_r}(\mathsf{h})$ can be written in the power form  $(\Psi^{(1)}_{1,\rho_l,\rho_r}(\mathsf{h}))^N$, and more precisely, \eqref{t1} tells us that
\begin{eqnarray*}
\Psi^{(1)}_{N,\rho_l,\rho_r}(\mathsf{h}) &=& \left(\frac{1}{(\rho_r-\rho_l)(1-e^{\mathsf{h}})} \log \frac{1+(1-
e^{\mathsf{h}})\rho_r}{1+(1-e^{\mathsf{h}})\rho_l}\right)^N \label{unmezzo}
\end{eqnarray*}
As a consequence
\be
\frac 1 N \log \Psi^{(1)}_{N,\rho_l,\rho_r}(\mathsf{h}) = \log\left(\frac{1}{(\rho_r-\rho_l)(1-e^{\mathsf{h}})} \log \frac{1+(1-e^{\mathsf{h}})\rho_r}{1+(1-e^{\mathsf{h}})\rho_l}\right) \;.
\ee
In other words, for  $s=1/2$ the finite volume pressure does not depend on $N$ and it coincides with the pressure at infinite volume.

\er

\subsection{Equivalence between additivity principle and variational problem}
\label{section-constant2}

In this section we will prove that  the fact that the modified pressure
\be\label{tPab}
\tilde{P}^{[a,b]}_{\rho_a,\rho_b} (h) := P^{[a,b]}_{\rho_a,\rho_b} (h) +
2s (b-a) \log \left(\frac{\rho_b-\rho_a}{b-a}\right)
\ee
 satisfies the additivity principle \eqref{adgeneral}, combined with the continuity of $\tilde{P}^{[a,b]}_{\rho_a,\rho_b}$ with respect to convergence in $L^1$ and with formula \eqref{pressab}
 that gives an explicit expression of the action of $P^{[a,b]}_{\rho_a,\rho_b} (h)$ on constant functions $h(x)=\mathsf{h}$ for all $x\in [a,b]$, allows to identify  the pressure functional $P^{[0,1]}_{\rho_l,\rho_r}$ on a generic function $h:[0,1]\to 1$, $h\in C^1$ as the solution of the variational problem:
 \begin{align}
P^{[0,1]}_{\rho_l,\rho_r} (h)
&=\tilde{P}^{[0,1]}_{\rho_l,\rho_r} (h)-
2s  \log\left({\rho_r-\rho_l}\right)\nn\\
&=2s \cdot \sup_{\theta}
\int_0^1  dx \,\left[ \log \left( \frac {\theta'(x)}{\rho_r-\rho_l} \right)+\log \left( \frac 1{1+\theta(x)(1-e^{h(x)})} \right)\right] \;.
  \end{align}

\subsubsection{Variational problem implies additivity  principle}
Consider $0=x_0< x_1 < \ldots < x_\kappa =1$.
Assume $h(x) = \sum_{i=1}^\kappa {h}_i(x)\indic{[x_{i-1},x_i]}(x)$ for $x\in[0,1]$
where ${h}_i$ is the restriction of $h$ to the interval $[x_{i-1},x_i]$.
Then the MFT variational problem
can be written as follows:
$$
 P^{[0,1]}_{\rho_l,\rho_r} (h) = \sup_{\theta}
\sum_{i=1}^\kappa \int_{x_{i-1}}^{x_i}
2s \left[
 \log \Big(
 \frac{1}{1+(1-e^{{h}_i(x)})\theta(x)} \Big)
+  \log\Big(
 \frac{\theta'(x)}{(\rho_r -\rho_l)}\Big) \right]  dx
$$
where the supremum is over monotonic $C^1$ functions $\theta:[0,1]\to \mathbb{R}$  such that $\theta(0)=
\rho_l$ and $\theta(1)=  \rho_r$.
Equivalently we can write
$$
 P^{[0,1]}_{\rho_l,\rho_r} (h) =
 \sup_{\rho_l= \rho_0 < \rho_1 < \ldots <\rho_\kappa =\rho_r}
\sum_{i=1}^\kappa \sup_{\theta_i} \int_{x_{i-1}}^{x_i} 2s
 \left[
 \log \Big(
 \frac{1}{1+(1-e^{{h}_i(x)})\theta_i(x)} \Big)
+  \log\Big(
 \frac{\theta_i'(x)}{\rho_r -\rho_l}\Big) \right] dx
$$
where the $i^{\text{th}}$ supremum is now over monotone $C^1$ functions $\theta_i:[x_{i-1},x_i]\to \mathbb{R}$  such that $\theta_i(x_{i-1})=
\rho_{i-1}$ and  $\theta_i(x_i)= \rho_i$.

We now write the right hand side above in terms of the pressures of each interval, i.e.
\begin{align}
\label{828}
 P^{[0,1]}_{\rho_l,\rho_r} (h) =
 \sup_{\rho_l= \rho_0 < \rho_1 < \ldots <\rho_\kappa =\rho_r}
\sum_{i=1}^\kappa
&  \Big[
\sup_{\theta_i} \int_{x_{i-1}}^{x_i} 2s
 \left[
 \log \Big(
 \frac{1}{1+(1-e^{{h}_i(x)})\theta_i(x)} \Big)
+  \log\Big(
 \frac{(x_i-x_{i-1})\theta_i'(x)}{\rho_i -\rho_{i-1}}\Big)
 \right] dx
\nn\\
&
+ 2s (x_{i}-x_{i-1}) \log\Big(
 \frac{\rho_i-\rho_{i-1}}{(x_i-x_{i-1}) (\rho_r -\rho_l)}\Big)
 \Big].
\end{align}
Define the pressure of the volume $[a,b]$ with boundary parameters $\rho_a,\rho_b$ as
$$
P^{[a,b]}_{\rho_a,\rho_b}(h) = \sup_{\theta}
 \int_{a}^{b} 2s
 \left[
 \log \Big(
 \frac{1}{1+(1-e^{{h}(x)})\theta(x)} \Big)
+  \log\Big(
 \frac{(b-a)\theta'(x)}{\rho_b -\rho_{a}}\Big)
 \right] dx
$$
where the supremum is  over monotone $C^1$ functions $\theta:[a,b]\to \mathbb{R}$  such that $\theta(a)=
\rho_a$ and  $\theta(b)= \rho_b$. Then \eqref{828} can be written as
\begin{align}
P^{[0,1]}_{\rho_l,\rho_r} (h)
& = \sup_{\rho_l= \rho_0 < \rho_1 < \ldots <\rho_\kappa =\rho_r}
\sum_{i=1}^\kappa \Big[P^{[x_{i-1},x_i]}_{\rho_{i-1},\rho_i} ({h}_i)  + 2s (x_{i}-x_{i-1}) \log\Big(
 \frac{\rho_i-\rho_{i-1}}{(x_i-x_{i-1}) (\rho_r -\rho_l)}\Big)\Big] \;. \nn
\end{align}
As a consequence,
we obtain that the modified pressure \eqref{tPab} fulfills the additivity principle
\begin{align}
\tilde P^{[0,1]}_{\rho_l,\rho_r} (h)
& = \sup_{\rho_l= \rho_0 < \rho_1 < \ldots <\rho_\kappa =\rho_r}
\sum_{i=1}^n\tilde P^{[x_{i-1},x_i]}_{\rho_{i-1},\rho_i} (h_i) \ .
\end{align}

 \subsubsection{Additivity principle  implies variational problem}

For any $C^1$ function $h:[0,1]\to \mathbb R$ we can produce a discretization by fixing a sequence of piecewise constant functions $h^{(\kappa)}:[0,1]\to \mathbb R$ defined as follows:
\be
h^{(\kappa)}(x)=\sum_{i=1}^{\kappa} \mathsf{h}_i \cdot \indic{[x_{i-1}, x_i]}(x),\qquad \kappa\in \mathbb N, \quad \mathsf{h}_1, \ldots, \mathsf{h}_{\kappa} \in \mathbb R
\ee
where
\be \label{xi}
x_i=\tfrac i \kappa \qquad \text{and}\quad \mathsf{h}_i:= h(x_i)=h\left(\tfrac i\kappa\right)
\ee
so that
\be
h^{(\kappa)}(x)= h\left(\tfrac {\lceil \kappa x\rceil}\kappa\right).
\ee
Then we have that $h^{(\kappa)}$ converges to $h$ in $L^1$.
We can define an analogous approximation for any $C^1$ function $\theta:[0,1]\to \mathbb R$ that is non-decreasing and such that $\theta (0)=\rho_l$ and $\theta(1)=\rho_r$.
We do it by defining the piecewise constant functions
\be
\theta^{(\kappa)}(x)=\sum_{i=1}^{\kappa} \rho_i \cdot \indic{[x_{i-1}, x_i]}(x),\qquad \text{for} \quad \rho_i:= \theta(x_i)=\theta\left( \tfrac i\kappa\right)
\ee
so that $\rho_l= \rho_0 < \rho_1 < \ldots <\rho_\kappa =\rho_r$ and
\be
\theta^{(\kappa)}(x)= \theta\left(\tfrac {\lceil \kappa x\rceil}\kappa\right).
\ee
We  assume that the modified pressure \eqref{tPab}  satisfies  the additivity principle \eqref{adgeneral} and apply this property to the case in which the external field is the piecewise constant function $h^{(\kappa)}$:

\begin{align}\label{addp}
\tilde P^{[0,1]}_{\rho_l,\rho_r} (h^{(\kappa)})
& = \sup_{\rho_l= \rho_0 < \rho_1 < \ldots <\rho_\kappa =\rho_r}
\sum_{i=1}^\kappa\tilde P^{[x_{i-1},x_i]}_{\rho_{i-1},\rho_i} (\mathsf{h}_i)
\end{align}
where
\be
\tilde{P}^{[x_{i-1},x_i]}_{\rho_{i-1},\rho_i} (\mathsf{h}_i) := P^{[x_{i-1},x_i]}_{\rho_{i-1},\rho_i} (\mathsf{h}_i) +
2s (x_i-x_{i-1}) \log \left(\frac{\rho_i-\rho_{i-1}}{x_i-x_{i-1}}\right).
\ee
We can use now  formula \eqref{pressab} which gives the pressure functional on constant functions
\be
\label{pressi}
P^{[x_{i-1},x_i]}_{\rho_{i-1},\rho_i} (\mathsf{h}_i) =2s(x_i-x_{i-1})\log \left( \frac 1{(\rho_i-\rho_{i-1})\left(1-e^{\mathsf{h}_i}\right)} \cdot  \log\frac{1+\rho_i(1-e^{\mathsf{h}_i})}{1+\rho_{i-1}(1-e^{\mathsf{h}_i})}\right)\nn
\ee
from which we compute
\be
\tilde{P}^{[x_{i-1},x_i]}_{\rho_{i-1},\rho_i} (\mathsf{h}_i) :=2s (x_i-x_{i-1}) \log \left( \frac 1{(x_i-x_{i-1})\left(1-e^{\mathsf{h}_i}\right)} \cdot  \log\frac{1+\rho_i(1-e^{\mathsf{h}_i})}{1+\rho_{i-1}(1-e^{\mathsf{h}_i})}\right).
\ee
Using \eqref{addp} and \eqref{xi} we have
\begin{align}
\tilde P^{[0,1]}_{\rho_l,\rho_r} (h^{(\kappa)})
& =  \sup_{\rho_l= \rho_0 < \rho_1 < \ldots 	<\rho_\kappa=\rho_r}
2s\sum_{i=1}^\kappa (x_i-x_{i-1})  \log \left( \frac 1{(x_i-x_{i-1})\left(1-e^{\mathsf{h}_i}\right)} \cdot  \log\frac{1+\rho_i(1-e^{\mathsf{h}_i})}{1+\rho_{i-1}(1-e^{\mathsf{h}_i})}\right)\nn \\
& =\sup_{\rho_l= \rho_0 < \rho_1 < \ldots <\rho_\kappa =\rho_r} \frac{2s}\kappa
\sum_{i=1}^\kappa \log \left( \frac \kappa{\left(1-e^{\mathsf{h}_i}\right)} \cdot  \log\frac{1+\rho_i(1-e^{\mathsf{h}_i})}{1+\rho_{i-1}(1-e^{\mathsf{h}_i})}\right)\label{quoo}.
\end{align}
Writing
\be
\log\left(\frac{1+\rho_i(1-e^{\mathsf{h}_i})}{1+\rho_{i-1}(1-e^{\mathsf{h}_i})}\right)= \log \left(1+\frac{(\rho_i-\rho_{i-1})(1-e^{\mathsf{h}_i})}{1+\rho_{i-1}(1-e^{\mathsf{h}_i})}\right)
\ee
and approxamiting
\be
\rho_i-\rho_{i-1}= \theta \left(\tfrac{\lceil \kappa x\rceil}{\kappa}\right)-\theta \left(\tfrac{\lceil \kappa x\rceil-1}{\kappa}\right)=\tfrac 1 \kappa \,\theta'(x) + o\left(\tfrac 1\kappa \right) \qquad \text{for} \qquad x_{i-1}\le x<x_i
\ee
and
\be
\mathsf{h}_i=h\left(\tfrac{\lceil \kappa x\rceil}{\kappa}\right) =h(x) + o(1) \qquad \text{for} \qquad x_{i-1}\le x<x_i
\ee
we get 
\be
\frac{(\rho_i-\rho_{i-1})(1-e^{\mathsf{h}_i})}{1+\rho_{i-1}(1-e^{\mathsf{h}_i})}=\frac 1 \kappa \cdot \frac{\theta'(x)(1-e^{h(x)})}{1+\theta(x)(1-e^{h(x)})} + o\left(\frac 1\kappa \right) \qquad \text{for} \qquad x_{i-1}\le x<x_i
\ee
and, as a consequence, taking the Taylor expansion of $\log(1+x)$ we obtain
\be
\log \left(1+\frac{(\rho_i-\rho_{i-1})(1-e^{\mathsf{h}_i})}{1+\rho_{i-1}(1-e^{\mathsf{h}_i})}\right)=\frac 1 \kappa \cdot \frac{\theta'(x)(1-e^{h(x)})}{1+\theta(x)(1-e^{h(x)})} + o\left(\frac 1\kappa \right) \qquad \text{for} \qquad x_{i-1}\le x<x_i \;.
\ee
Substituting this in \eqref{quoo} and taking the limit as $\kappa\to \infty$, via convergence of the Riemann sum to the corresponding integral we obtain that
\begin{align}
\tilde P^{[0,1]}_{\rho_l,\rho_r} (h) =\lim_{\kappa\to \infty}\tilde P^{[0,1]}_{\rho_l,\rho_r} (h^{(\kappa)})
& =2s \cdot \sup_{\theta}
\int_0^1 \log \left( \frac{\theta'(x)}{1+\theta(x)(1-e^{h(x)})} \right)
\end{align}
where the first identity follows from the continuity of the modified pressure functional with respect to convergence of function in $L^1$.
Now, using again \eqref{tPab} we conclude that
\begin{align}
P^{[0,1]}_{\rho_l,\rho_r} (h)
&=\tilde{P}^{[0,1]}_{\rho_l,\rho_r} (h)-
2s  \log\left({\rho_r-\rho_l}\right)\nn\\
&=2s \cdot \sup_{\theta}
\int_0^1 dx \, \left[ \log \left( \frac {\theta'(x)}{\rho_r-\rho_l} \right)+\log \left( \frac 1{1+\theta(x)(1-e^{h(x)})} \right)\right].
  \end{align}

\subsection{Finite volume}
\label{section-constant3}

In what follows we show that the moment generating function $\Psi_{N,\rho_l,\rho_r}$  has another expression which differs from the ones in terms of $N$-fold sums and $N$-folds integrals of Sections \ref{Nfoldsumssubsection} and \ref{N-fold integrals}. To some extent this expression is more clear because it only relies on finite sums.

\subsubsection{Recurrence relation}
We start from
the integral equation \eqref{integral-eq333}
relating partition
functions of different sizes and specialise it to the case $N_1=1$ and $N_2=N-1$. This becomes
\begin{eqnarray}
\label{recu-psi}
\Psi_{N,\rho_l,\rho_r}( \mathsf{h}_1,\ldots,  \mathsf{h}_N)
&= &
\int_{0}^1 du \Big( \frac{1}{1+(1-e^{ \mathsf{h}_1})\theta(u)}\Big)^{2s}
\Psi_{N-1,\theta(u),\rho_r}( \mathsf{h}_2,\ldots,  \mathsf{h}_{N})\nn\\
& & \qquad
\frac{\Gamma(2s(N+1))}{\Gamma(2s)\Gamma(2s N)} u^{2s-1} (1-u)^{2sN-1} \;.
\end{eqnarray}
Thanks to the relation \eqref{Phi-fct} between $\Psi_{N,\rho_l,\rho_r}$ and $\Phi_N$ we can turn \eqref{recu-psi} in a recurrence relation for the function $\Phi_N$, namely
\begin{eqnarray}
\Phi_{N}(c_1,\ldots,c_N)
&= &
\int_{0}^1 du \Big( \frac{1}{1-(1-u)c_1}\Big)^{2s}
\Phi_{N-1}\Big((1-u)c_2,\ldots,(1-u) c_{N}\Big)\nn\\
& & \qquad
\frac{\Gamma(2s(N+1))}{\Gamma(2s N)\Gamma(2s)} u^{2s-1} (1-u)^{2sN-1} \;.
\end{eqnarray}
Changing the integration variable to $t = 1-u$ one obtains
\begin{equation}\label{osso}
\Phi_{N}(c_1,\ldots,c_N) =  \frac{1}{B(2sN,2s)}  \int_{0}^1 dt \left(\frac{1}{1-tc_1}\right)^{2s} t^{2sN-1} (1-t)^{2s -1} \Phi_{N-1}(t c_2,\ldots,t c_N)
\end{equation}
where $B(2sN, 2s) = \frac{\Gamma(2sN) \Gamma(2s)}{\Gamma(2s(N+1))}$ is the Beta function.
\vskip.1cm
\noindent
Choosing a  constant external field $( \mathsf{h}, \ldots,  \mathsf{h})$ corresponds to choosing a vector $c_{N,\rho_l,\rho_r}(\mathsf{h}) $ (see \eqref{ci-fct}) with constant components $c_i := \left(  c_{N,\rho_l,\rho_r}\left( \mathsf{h}\right) \right)_{i}  =c \in \mathbb R$ for $i = 1, \ldots, N$. For convenience we use the notation $\Phi_N^{(1)}$ for the function:
\be
\Phi_N^{(1)}(c):= \Phi_N(c, \ldots, c)
\ee
then, specialising \eqref{osso} to the case $c_1= \ldots =  c_N=c \in \mathbb R$ we deduce the following recurrence relation on $\Phi_N^{(1)}$
\begin{align}
\label{recurr}
 \Phi^{(1)}_{N}(c) & = \frac{1}{B(2sN,2s)} \int_{0}^1 dt \left(\frac{1}{1-ct}\right)^{2s} t^{2sN-1} (1-t)^{2s -1} \Phi^{(1)}_{N-1}(t c) \\ &
= \mathbb{E}\left[ \left(  \dfrac{1}{1-c\mathfrak{B}}\right)^{2s} \Phi^{(1)}_{N-1}(c\mathfrak{B})    \right] \nn
\end{align}
where  the random variable $\mathfrak{B}$ is distributed as a Beta$(2sN, \ 2s)$. 
Now we will see that it is possible to turn the integral in the right hand side of \eqref{recurr} into a convolution.
To this aim, we perform the following change of variables $c = 1- e^{-2v}$ and define the random variable $\mathfrak{Z}$ via the relation $c \mathfrak{B} = 1-e^{-2 \mathfrak{Z}} $. Then  the density function of $\mathfrak{Z}$ is
\begin{equation}
f_{\mathfrak{Z}}(z) = \dfrac{1}{B(2sN, 2s)} \left( \dfrac{1-e^{-2z}}{1-e^{-2v}}\right) ^{2sN -1} \left( \dfrac{e^{-2z}-e^{-2v}}{1-e^{-2v}}\right) ^{2s -1} \dfrac{2e^{-2z}}{1-e^{-2v}} \;,
\end{equation}
which allows to rewrite the recurrence relation in \eqref{recurr}  as
\begin{align} \label{recurrphi}
 \Phi_{N}^{(1)}(1- e^{-2v}) & = \mathbb{E}\left[ \left( e^{\mathfrak{2Z}} \right)^{2s}
 \Phi^{(1)}_{N-1}(1- e^{-2\mathfrak{Z}})   \right] \;.
\end{align}
Using the density function of $\mathfrak{Z}$, the expression above can be conveniently rewritten as
\begin{align}
&\dfrac{e^{v(2s-1)} (1-e^{-2v})^{2s(N+1) -1} }{2^{2sN}}\cdot \Phi^{(1)}_{N}(1- e^{-2v}) =\nn\\
 &\dfrac{1}{B(2sN, 2s)} \int_{0}^{v} dz \ \dfrac{e^{z(2s-1)} (1-e^{-2z})^{2sN -1}}{2^{2s(N-1)}} \cdot  \Phi^{(1)}_{N-1}(1- e^{-2z}) \left(  \sinh(v-z) \right) ^{2s-1} \;.
\end{align}
Defining the l.h.s. above as
\begin{align} \label{G-Phi}
\mathrm{G}_N(v) := \dfrac{e^{v(2s-1)} (1-e^{-2v})^{2s(N+1) -1} }{2^{2sN}}\cdot \Phi^{(1)}_{N}(1- e^{-2v})
\end{align}
allows to read the recurrence relation as a convolution, i.e.
\begin{align}
\mathrm{G}_N(v) = \dfrac{1}{B(2sN, 2s)} \int_{0}^{v} dz \ \mathrm{G}_{N-1}(z) \left(  \sinh(v-z) \right) ^{2s-1} \;,
\end{align}
with $G_0 (v) = 2^{2s-1} \left(  \sinh(v) \right) ^{2s-1} $.
Iterating $N+1$ times, we can write $\mathrm{G}_N$ as
\begin{equation} \label{convolution}
\mathrm{G}_N(v) = 2^{2s-1} \dfrac{\Gamma(2s(N+1))}{\Gamma(2s)^{N+1}}  \left( g \ast \ldots \ast g \right)  (v)
\end{equation}
where the convolution is taken $N+1$ times and $g(v)=  \left(  \sinh(v) \right) ^{2s-1} $.
\subsubsection{Pressure via inverse Laplace transform}
In the previous section we have seen how the recurrence relation for the function $\Phi_{N}^{(1)}$ in \eqref{recurrphi} simplifies in a convolution relation for the function $\mathrm{G}_N $ 
in \eqref{convolution}.
Therefore, we denote by $\widehat{g} (\alpha) = \Laplace\{ g(v)\} (\alpha) $ the Laplace transform of the function $g(v)$ so that, when considering the Laplace transform on both sides of \eqref{convolution}, we get
\begin{equation} \label{laplace-convolution}
\widehat{\mathrm{G}}_N (\alpha)  =  2^{2s-1} \dfrac{\Gamma(2s(N+1))}{\Gamma(2s)^{N+1}} \left( \widehat{g} (\alpha)  \right) ^{N+1} \;.
\end{equation}
Computing the Laplace transform of $g(v)$, allows to explicitly write
\begin{equation} \label{invlapG}
\widehat{\mathrm{G}}_N (\alpha) 
=\ 2^{2s-1} \dfrac{\Gamma(2s(N+1))}{2^{2s(N+1)}} \left( \dfrac{\Gamma\left(  \frac{ \alpha + 1 -2s}{2} \right)}{\Gamma\left(  \frac{ \alpha + 1 +2s}{2} \right)} \right) ^{N+1}  \;.
\end{equation}
At this point it is clear that anti transforming $\widehat{\mathrm{G}}_N (\alpha) $ and using equation \eqref{G-Phi}, one can explicitly get an expression for the finite volume pressure for all $N$. This is the content of the next proposition.
\bp[Closed formula]
For $2s \in \mathbb{N}$, a closed formula for $\Phi^{(1)}_{N}$ given in terms of a finite sum  is
\begin{align} \label{finite}
\Phi^{(1)}_{N}(c) =\left(  \dfrac{2}{c} \right)^{2s(N+1)-1}\Gamma(2sN+2s) \sum_{j=0}^{2s-1} \sum_{k=0}^{N}
\dfrac{(-\log(1-c))^{N-k} }{2^{N-k} (N-k)! k!} \left( 1-c \right)^{j}
 \phi_{j,k}(\alpha_j)
\end{align}
where
\be
\phi_{j,k}(\alpha_j) = \sum_{j_0 + j_1 + \ldots + j_{2s-1} = k} {{k} \choose{j_0, j_1, \ldots , j_{2s-1}}} \prod_{\substack{i = 0 \\ i \neq j}}^{2s - 1} (-1)^{j_i} \dfrac{(N+j_i)!}{N!} (2i - 2j)^{-N-j_i-1} \;.
\ee
\ep
\begin{proof}
First we compute the inverse Laplace transform of $\widehat{\mathrm{G}}_N (\alpha) $ then we use  \eqref{G-Phi} to get the expression above. In order to invert the Laplace transform we notice that
\begin{equation*}
\dfrac{\Gamma\left(\frac{\alpha -2s +1}{2} \right) }{\Gamma\left(\frac{\alpha +2s +1}{2} \right)} = \prod_{i=0}^{2s-1} \dfrac{2^{2s}}{ \left( \alpha - (2s-1) +2i \right) }
\end{equation*}
in other words, $\widehat{\mathrm{G}}_N (\alpha)$ has  $2s$ poles, all with multiplicity $N+1$ namely $\alpha_i = 2s-1-2i $ for $i=0,\ldots, 2s-1$. The inverse Laplace transform of a rational function can be computed (see for example formula (21) of \cite{erdelyi}); in our case
\begin{equation*}
\widehat{\mathrm{G}}_N (\alpha)  =   2^{2s-1}\Gamma(2s(N+1)) \prod_{i=0}^{2s-1} \left( \dfrac{1}{\alpha - (2s-1) +2i } \right) ^{N+1} \;
\end{equation*}
has inverse Laplace transform
\begin{equation*}
\mathrm{G}_N(v) = 2^{2s-1}  \Gamma(2sN+2s) \sum_{j=0}^{2s-1}\sum_{k=0}^{N}  \dfrac{v^{N-k}}{(N-k)! k! }  \phi_{j,k}(\alpha_{j}) e^{\alpha_{j}v} \qquad \text{where} \qquad   \phi_{j,k} (\alpha) =  \dfrac{\partial^{k}}{\partial \alpha^{k}}   \prod_{\substack{i = 0 \\ i \neq j}}^{2s - 1} \left( \dfrac{1}{\alpha - \alpha_i} \right)^{N+1} \;.
\end{equation*}
Now we show an explicit formula for the factors $\phi_{j,k}(\alpha_{j}) $, which can be computed using the general Leibniz rule for the product of functions, i.e.
\begin{equation*}
\phi_{j,k}(\alpha) = \sum_{j_0 + j_1 + \ldots + j_{2s-1} = k} {{k} \choose{j_0, j_1, \ldots , j_{2s-1}}} \prod_{\substack{i = 0 \\ i \neq j}}^{2s - 1} \dfrac{\partial^{j_i} }{\partial \alpha^{j_i}} \left(\dfrac{1}{\alpha - \alpha_i} \right)^{N+1}
\end{equation*}
where ${{k} \choose{j_0, j_1, \ldots , j_{2s-1}}} $ is the multinomial coefficient
 and the ${j_i}^{th}$ derivative with respect to $\alpha$ is
\begin{equation*}
\dfrac{\partial^{j_i} }{\partial \alpha^{j_i}} \left(\dfrac{1}{\alpha - \alpha_i} \right)^{N+1} = (-1)^{j_i} \dfrac{(N+j_i)!}{N!} (\alpha - \alpha_i)^{-N-1-j_i}
\end{equation*}
so that
\be
\phi_{j,k}(\alpha) = \sum_{j_0 + j_1 + \ldots + j_{2s-1} = k} {{k} \choose{j_0, j_1, \ldots , j_{2s-1}}} \prod_{\substack{i = 0 \\ i \neq j}}^{2s - 1} (-1)^{j_i} \dfrac{(N+j_i)!}{N!} (\alpha - \alpha_i)^{-N-j_i-1} \;.
\ee
All in all, recalling that the residues are $\alpha_j = 2s-1-2j$  we get an explicit expression for $\mathrm{G}_N(v)$,
\begin{align*}
\mathrm{G}_N(v)= 2^{2s-1}  \Gamma(2sN+2s) \sum_{j=0}^{2s-1}\sum_{k=0}^{N}  \dfrac{ v^{N-k} \ e^{(2s-2j-1)v }}{(N-k)! k! } \phi_{j,k}(\alpha_j) \;.
\end{align*}
The expression in equation \eqref{finite} is then obtained from \eqref{G-Phi} setting $c = 1-e^{-2v}$ and rewriting for $ \Phi^{(1)}_{N} (c) $,
\begin{equation}\label{phiexplicit}
\Phi^{(1)}_{N} (c) = \dfrac{2^{2sN}}{c^{2s(N+1) -1}} (1-c)^{s-1/2} \ \mathrm{G}_N \left(-\frac{1}{2} \log\left( 1-c \right)  \right) \;.
\end{equation}
\end{proof}

We now show how the above computation for the moment generating function specialises for the first two cases $s=1/2$ and $s=1$.
\paragraph{Case $s=1/2$.}
In this first case one can check that equation \eqref{invlapG} simplifies to
\begin{equation*}
\widehat{\mathrm{G}}_N (\alpha) = \dfrac{N!}{\alpha^{N+1}}
\end{equation*}
and its inverse Laplace transform is
\begin{equation*}
\mathrm{G}_N (v) = v^{N} \;.
\end{equation*}
Recalling equation \eqref{G-Phi} and considering again the change of variable $c=1-e^{-2v}$, we obtain
\begin{equation} \label{thisone}
 \Phi^{(1)}_{N} (c) = \left(- \dfrac{1}{c} \log(1-c) \right)^{N}  \;.
\end{equation}
Notice that, using the map in \eqref{ci-fcti} with $s=1/2$ and constant $\mathsf{h}$ we recover the expression \eqref{unmezzo} for the MGF on a constant field $\Psi^{(1)}_{N, \rho_l,\rho_r}(\mathsf{h})$.

\paragraph{Case $s=1$.}
This is the first non-trivial value of $s$, notice that for all $s \neq 1/2 $  the Laplace transform we compute depends on exponential functions and computations are more involved. We proceed as before. From \eqref{invlapG} we can write
\begin{equation*}
\widehat{\mathrm{G}}_N (\alpha)  =  2 \  \Gamma (2N+2) \left( \dfrac{1}{\alpha^{2}-1} \right)^{N+1} \;.
\end{equation*}
The poles of $\widehat{\mathrm{G}}_N$ are $\alpha_0 = -1$ and $\alpha_1 = 1$, while its Laplace inverse is
\begin{equation*}
\mathrm{G}_N (v)  =  2 \  \Gamma (2N+2) \sum_{k=0}^{N}  \dfrac{v^{N-k}}{k! (N-k)!} \left[ \phi_{0,k}(-1) e^{- v }+  \phi_{1,k}(1) e^{v}\right]
\end{equation*}
where
\begin{align*}
\phi_{0,k}(-1) &= (-1)^{k}\dfrac{(N+k)!}{N!} (-2)^{-N-1-k} \qquad \text{and}  \qquad
\phi_{1,k}(1) = (-1)^{k}\dfrac{(N+k)!}{N!} (2)^{-N-1-k} \;.
\end{align*}
This leads to
\begin{equation*}
\mathrm{G}_N (v)  =  2 \  \Gamma (2N+2) \sum_{k=0}^{N}  v^{N-k} \dfrac{(N+k)!}{N! k! (N-k)!} (-1)^{k} 2^{-N-1-k}\left[ e^{v }+  e^{- v} (-1)^{N+1+k} \right]  
\end{equation*}
and using equation \eqref{phiexplicit} we obtain
\begin{align*}
 \Phi^{(1)}_{N}(c) = \dfrac{\Gamma(2N+2)}{c^{2N+1}} \sum_{k=0}^{N} (-1)^N  \dfrac{(N+k)!}{N!k!(N-k)!} (\log (1-c))^{N-k} \left[ 1+ (1-c)  (-1)^{N+k+1} \right]  \;.
  \end{align*}

\appendix
\section{Comparison and integral representation of moments}\label{app:compare}
In this appendix we show that \eqref{eq:measureee} coincides with the integral representation of the steady state in \eqref{mixtresult-explicit}. To do so, we consider the integral representation of the factorial moments that follows immediately from inserting \eqref{mixtresult-explicit} into \eqref{eq:Gasmu}. It reads
\begin{equation}\label{mixtresult-explicit2}
 \begin{split}
  G(\xi)
&=  \mathcal{N}(N,s)\cdot
 \int_{\rho_l}^{\rho_r}d\theta_1 \int_{\theta_1}^{\rho_r}d\theta_2\cdots  \int_{\theta_{N-1}}^{\rho_r}d\theta_N      \Big[\prod_{i=1}^{N+1} (\theta_{i}-\theta_{i-1})^{2s-1}\Big]
\Big[\prod_{i=1}^{N} \theta_i^{\xi_i}\Big] \,,
 \end{split}
\end{equation}
with the normalisation
\begin{equation}
 \mathcal{N}(N,s)=
\frac{\Gamma(2s(N+1))}{\Gamma(2s)^{N+1}}(\rho_r-\rho_l)^{1-2s(N+1)}\,.
\end{equation}
To show that \eqref{mixtresult-explicit2} coincides with the factorial moments of Theorem~\ref{SFM-new}, we consider the auxiliary function
\begin{equation}
 \mu''(\eta) =\sum_{\xi\geq\eta}\Big[\prod_{i=1}^N\frac{(-\rho_r)^{\eta_i-\xi_i}}{\eta_i!}  \binom{\eta_i}{\xi_i}   \frac{\Gamma(2s+\eta_i)}{\Gamma(2s)} \Big]G(\xi)\,,
\end{equation}
as introduced in \cite{fg}. It is written in terms of the integrals as
\begin{equation}\label{mixtresult-explicit3}
 \begin{split}
  \mu''(\eta)
&=
\mathcal{N}(N,s)\Big[\prod_{i=1}^{N}\frac{1}{\eta_i!}\frac{\Gamma(2s+\eta_i)}{\Gamma(2s)} \Big]
 \int_{\rho_l}^{\rho_r}d\theta_1 \int_{\theta_1}^{\rho_r}d\theta_2\cdots  \int_{\theta_{N-1}}^{\rho_r}d\theta_N
  \Big[\prod_{i=1}^{N+1} (\theta_{i}-\theta_{i-1})^{2s-1}\Big]
\Big[\prod_{i=1}^{N} \left(\theta_i-\rho_r\right)^{\eta_i}\Big]\,.
 \end{split}
\end{equation}
These integrals can be evaluated explicitly for any spin $s$ and length $N$. Introducing  the variables $u_i=\theta_i-\rho_r$ and $u_0=\rho_l-\rho_r$
and using repeatedly the formula
\begin{equation}
 \int_0^xdy\, y^a (y-x)^b=(-1)^b x^{a+b+1}\frac{ \Gamma (a+1) \Gamma (b+1)}{\Gamma (a+b+2)}\,.
\end{equation}
We find
\begin{equation}\label{mixtresult-explicit5}
 \begin{split}
  \mu''(\eta)
&= (\rho_l-\rho_r)^{|\eta|}\Big[\prod_{i=1}^{N}\frac{1}{\eta_i!}\frac{\Gamma(2s+\eta_i)}{\Gamma(2s)} \Big]
\frac{\Gamma(2s(N+1))}{\Gamma(2s(N+1)+|\eta|)}\cdot\prod_{k=1}^N  \frac{\Gamma(2s(N-k+1)+\sum_{i=k}^N\eta_{i})}{\Gamma(2s(N-k+1)+\sum_{i=k+1}^N\eta_{i})}\,,
 \end{split}
\end{equation}
which coincides with \cite[(6.3)]{fg}.

\end{document}